\documentclass[11pt]{article}
\usepackage{tocbibind} 
\usepackage{amsmath,amsthm,amssymb,thmtools}
\usepackage{fancyhdr}
\usepackage[british]{babel}
\usepackage{geometry,mathtools}
\usepackage{enumitem}
\usepackage{algpseudocode}
\usepackage{dsfont}
\usepackage{centernot}
\usepackage{xstring}
\usepackage{colortbl}
\usepackage{pifont}

\usepackage[section]{algorithm}
\usepackage{graphicx}
\usepackage[font={footnotesize}]{caption}
\usepackage[usenames,dvipsnames,table]{xcolor}
\usepackage{pstricks,tikz}
\usetikzlibrary{patterns,calc,fadings,decorations.pathreplacing}
\usepackage[h]{esvect}
\usepackage[
bookmarksopen=true,
bookmarksopenlevel=1,
colorlinks=true,
linkcolor=darkblue,
linktoc=page,
citecolor=darkblue,
urlcolor=darkblue
]{hyperref}
\usepackage{bbm}

\usepackage{multirow}

\numberwithin{equation}{section}

\definecolor{darkblue}{rgb}{0,0,0.5}

\newdimen\margin
\def\textno#1&#2\par{
   \margin=\hsize
   \advance\margin by -4\parindent
          \setbox1=\hbox{\sl#1}
   \ifdim\wd1 < \margin
      $$\box1\eqno#2$$
   \else
      \bigbreak
      \hbox to \hsize{\indent$\vcenter{\advance\hsize by -3\parindent
      \it\noindent#1}\hfil#2$}
      \bigbreak
   \fi}

\usepackage{cleveref}[capitalise, compress, nameinlink, noabbrev]

\theoremstyle{plain}
\newtheorem{theorem}{Theorem}[section]
\crefname{theorem}{Theorem}{Theorems}

\newtheorem{prop}[theorem]{Proposition}
\crefname{prop}{Proposition}{Propositions}

\crefname{cor}{Corollary}{Corollaries}

\newtheorem{lemma}[theorem]{Lemma}
\crefname{lemma}{Lemma}{Lemmata}

\newtheorem{conj}[theorem]{Conjecture}
\crefname{conj}{Conjecture}{Conjectures}

\newtheorem{problem}[theorem]{Problem}
\crefname{problem}{Problem}{Problems}

\newtheorem{observation}[theorem]{Observation}
\crefname{observation}{Observation}{Observations}

\crefname{setup}{Setup}{Setups}

\newtheorem{fact}[theorem]{Fact}
\crefname{fact}{Fact}{Facts}

\crefname{remark}{Remark}{Remarks}

\crefname{example}{Example}{Examples}

\theoremstyle{definition}
\newtheorem{defin}[theorem]{Definition}
\crefname{defin}{Definition}{Definitions}

\crefname{construction}{Construction}{Constructions}

\crefname{question}{Question}{Questions}

\numberwithin{equation}{section}

\crefname{section}{Section}{Sections}
\crefname{appendix}{Appendix}{Appendix}

\crefname{figure}{Figure}{Figures}

\def\proof{\removelastskip\penalty55\medskip\noindent\begin{stepenv}\end{stepenv}{\bf Proof. }} 

\def\lateproof#1{\removelastskip\penalty55\medskip\noindent\begin{stepenv}\end{stepenv}{\bf Proof of #1. }} 

\def\noproof{{\unskip\nobreak\hfill\penalty50\hskip2em\hbox{}\nobreak\hfill%
       $\square$\parfillskip=0pt\finalhyphendemerits=0\par}\goodbreak}
\def\endproof{\noproof\bigskip}

\def\claimproof{\removelastskip\penalty55\medskip\noindent{\em Proof of claim: }}

\newcounter{stepenv}
\newenvironment{stepenv}[1][]{\refstepcounter{stepenv}}{}

\newcounter{step}[stepenv]

\newcounter{substep}[step]
\renewcommand{\thesubstep}{\thestep.\arabic{substep}}

\newcounter{claim}[stepenv]
\newenvironment{claim}[1][]{\refstepcounter{claim}\par\medskip\noindent%
        \textit{Claim~\theclaim. #1} \itshape\rmfamily}{\medskip}

\newcommand{\cA}{\mathcal{A}}

\newcommand{\cE}{\mathcal{E}}

\newcommand{\cH}{\mathcal{H}}

\newcommand{\cL}{\mathcal{L}}
\newcommand{\cM}{\mathcal{M}}

\newcommand{\cR}{\mathcal{R}}
\newcommand{\cS}{\mathcal{S}}
\newcommand{\cT}{\mathcal{T}}

\newcommand{\cX}{\mathcal{X}}

\newcommand{\bN}{\mathbb{N}}

\newcommand{\defn}{\emph}

\newcommand{\prob}[1]{\mathrm{\mathbb{P}}\left[#1\right]}

\newcommand{\expn}[1]{\mathrm{\mathbb{E}}\left[#1\right]}

\def\sm{\setminus}
\newcommand{\Set}[1]{\{#1\}}

\def\In{\subseteq}

\def\COMMENT#1{}
\def\TASK#1{}
\let\TASK=\footnote             

\begin{document}

\title{On Kotzig's conjecture in random graphs}

\author{
Stefan Glock \thanks{Fakultät für Informatik und Mathematik, Universität Passau, Germany.
\emph{Email}: \href{mailto:stefan.glock@uni-passau.de}{\tt stefan.glock@uni-passau.de}, \href{mailto:amedeo.sgueglia@uni-passau.de}{\tt amedeo.sgueglia@uni-passau.de}.
SG is funded by the Deutsche Forschungsgemeinschaft (DFG, German Research Foundation) – 542321564;
AS is funded by the Alexander von Humboldt Foundation.}
\and
Amedeo Sgueglia \footnotemark[1]
}

\date{}

\maketitle

\begin{abstract} 
In 1963, Anton Kotzig famously conjectured that $K_{n}$, the complete graph of order $n$, where $n$ is even, can be decomposed into $n-1$ perfect matchings such that every pair of these matchings forms a Hamilton cycle.
The problem is still wide open and here we consider a variant of it for the binomial random graph $G(n,p)$. We prove that, for every fixed $k$, there exists a constant $C=C(k)$ such that, when $p\ge \frac{C \log n}{n}$, with high probability, $G(n,p)$ contains $k$ edge-disjoint perfect matchings with the property that every pair of them forms a Hamilton cycle. In fact, our main result is a very precise counting result for~$K_n$. We show that, given any $k$ edge-disjoint perfect matchings $M_1,\dots,M_k$, the probability that a uniformly random perfect matching $M^*$ in $K_n$ has the property that $M^*\cup M_i$ forms a Hamilton cycle for each $i\in [k]$ is $\Theta_k(n^{-k/2})$. This is proved by building on a variety of methods, including a random process analysis, the absorption method, the entropy method and the switching method. 
The result on the binomial random graph follows from a slight strengthening of our counting result via the recent breakthroughs on the expectation threshold conjecture.
\end{abstract}

\section{Introduction}
\label{sec:intro}

A \defn{$1$-factor} or \defn{perfect matching} of a graph $G$ is a collection of pairwise vertex-disjoint edges of $G$ covering all the vertices of $G$.
A partition of the edge set of $G$ into $1$-factors is called a \emph{$1$-factorisation} of $G$. 
The following standard construction shows that if $n$ is even, then the complete graph of order $n$, denoted by $K_n$, has a $1$-factorisation. 
Place all but one of the vertices at the corners of a regular $(n-1)$-gon, with the remaining vertex at the center. 
Observe that, with this arrangement of vertices, choosing an edge $e$ from the center to a polygon vertex together with all possible edges that lie on lines perpendicular to $e$ gives a $1$-factor. The $1$-factors that can be constructed in this way form a $1$-factorisation of $K_n$.

\begin{figure}[htp]
    \centering
    \includegraphics[scale=1]{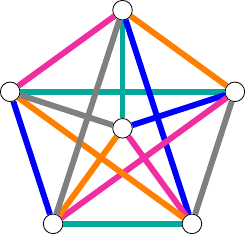}
    \caption{A perfect $1$-factorisation of $K_6$, where each $1$-factor is drawn with a different colour.}
    \label{fig:perfect_factorisation}
\end{figure}

Kotzig observed that if $n-1$ is a prime, then each pair of $1$-factors in the construction above induces a Hamilton cycle (cf.~\Cref{fig:perfect_factorisation}).
A $1$-factorisation in which any two distinct $1$-factors induce a Hamilton cycle is called a \defn{perfect $1$-factorisation}.
Motivated by this, in 1963, at ``the first international conference devoted to graph theory" (cf.~\cite{KO:14c}), Kotzig made the following conjecture, which has become known as the \defn{perfect $1$-factorisation conjecture}.
\begin{conj}[Kotzig's conjecture~\cite{kotzig:64}]
\label{conj:kotzig}
    Let $n \ge 2$ be an even integer. Then $K_n$ has a perfect $1$-factorisation.
\end{conj}

Kotzig's observation proves the conjecture for all $n$ such that $n-1$ is a prime number.
Later, Anderson~\cite{anderson:73} showed that the conjecture holds for all $n$ such that $n/2$ is a prime.
Much effort has been put into attacking this conjecture but it is still very far from being resolved. Historically, each increase in the smallest order for which Kotzig’s conjecture was unresolved had merited a paper.
In fact, besides the general results mentioned above, the conjecture is known to be true for some more sporadic values of $n$ and the current smallest open case is $n=64$.
We refer to the survey of Rosa~\cite{rosa:2019} for the full list of known values (with the case $n=54$ only appearing in a recent paper of Pike~\cite{pike:2019}).
In addition to the existence statement, there has been also some interest in counting the number of non-isomorphic perfect $1$-factorisations of $K_n$, often using computer-assisted proofs. 
The largest value of $n$ for which the exact answer is known is $n=16$: Gill and Wanless~\cite{GW:20} and, independently, Meszka~\cite{meszka:20} showed there are $3155$ non-isomorphic perfect $1$-factorisations of $K_{16}$.

\cref{conj:kotzig} has connections to several other problems.
A proper colouring of the edges of a graph $G$ is called \defn{acyclic} if there is no $2$-coloured cycle in $G$, that is, if the union
of any two colour classes induces a subgraph of $G$ which is a forest.
The \defn{acyclic edge chromatic number} of $G$, denoted by $a'(G)$, is the least number of colours in an acyclic edge colouring of $G$.
Fiam{\v c}{\'i}k~\cite{fiamcik:78} and, independently, Alon, Sudakov and Zaks~\cite{ASZ:01} conjectured that $a'(G) \le \Delta(G)+2$ for all graphs~$G$, where $\Delta(G)$ denotes the maximum degree of~$G$.
As shown in~\cite{ASZ:01}, if \cref{conj:kotzig} holds, then $a'(K_{2n})=a'(K_{2n+1})=2n+1$ for every $n$, which implies that the above conjecture is true for complete graphs.

Coming back to Kotzig's conjecture, Wanless~\cite{wanless:99} made a similar one for complete balanced bipartite graphs, stating that $K_{n,n}$ admits a perfect $1$-factorisation if and only if $n=2$ or $n$ is odd.
That the condition on $n$ is necessary may not look obvious at first sight, but it is well known (see the discussion at the end of \cref{sec:bipartite}) and goes back at least to Laufer~\cite{laufer:80}, who also showed that if $K_{n+1}$ has a perfect $1$-factorisation, then so does $K_{n,n}$.
Therefore, if true, \cref{conj:kotzig} implies its bipartite version.
One reason for interest in the bipartite analogue of \Cref{conj:kotzig} is its relationship with row-Hamiltonian Latin squares.
Let $m$ and $n$ be positive integers with $m \le n$. An \emph{$m \times n$ Latin rectangle} is an $m \times n$ matrix of $n$ symbols, each of which occurs exactly once in each row and at most once in each column. A \emph{Latin square} of order $n$ is an $n \times n$ Latin rectangle.
Suppose that $i$ and $j$ are two distinct rows of a Latin square $\cL$ with symbol set $S$ and define a permutation $\pi_{i,j} : S \to S$ by $\pi_{i,j}(\cL_{i,k}) = \cL_{j,k}$ for each $k \in S$. A Latin square of order $n$ is \emph{row-Hamiltonian} if $\pi_{i,j}$ consists of a single cycle for all~$i,j$.
It was noted~\cite{wanless:99} that a Latin square of order $n$ is row-Hamiltonian if and only if it contains no $a \times b$ Latin subrectangle with $1 < a \le b < n$.
This is an extremely rare property as it is known
(see, e.g.,~\cite{KSS:22}) that the overwhelming majority of Latin squares have quadratically many Latin subsquares.
There is a natural equivalence between row-Hamiltonian Latin squares of order $n$ and ordered perfect $1$-factorisations of $K_{n,n}$ (see~\cite{wanless:99}, where the argument is spelt out in detail), and thus the latter may be used to find explicit constructions of the former.
The conjecture of Wanless is a long way from being resolved. There are few known infinite families of perfect $1$-factorisations of $K_{n,n}$~\cite{AW:24,BMW:06,BMW:02}, and these only cover the case $n \in \{p, 2p - 1, p^2\}$ for some odd prime $p$.

Moreover, given the difficulty of proving \cref{conj:kotzig}, several notions that relax the perfectness requirement have appeared in the literature.
For example, resolving a conjecture of Archdeacon, it has been proved by Kr{\'a}lovi{\v c} and Kr{\'a}lovi{\v c}~\cite{KK:05} that for every even integer $n \ge 2$, there is a $1$-factorisation $\{M_1,\dots,M_{n-1}\}$ of $K_n$ such that $M_1 \cup M_j$ induces a Hamilton cycle for all $j \ge 2$. 
(These are known as \defn{semi-perfect $1$-factorisations}.)
Moreover, Wagner~\cite{wagner:92} studied an approximate version of \cref{conj:kotzig} in the following sense.
For a $1$-factorisation $\cM$ of $K_n$, let $c(\cM)$ be the number of pairs of $1$-factors of $\cM$ such that their union induces a Hamilton cycle.
Then define $c(n):=\max c(\cM)$, where the maximum is over all $1$-factorisations of $K_n$. 
Note that \cref{conj:kotzig} is equivalent to the statement that $c(n)=\binom{n-1}{2}$ for every even $n$, and Wagner offered some bounds on $c(n)$ in terms of the Euler's totient function.

Finally, one of the results we will prove says that, given a collection $\cM$ of perfect matchings on the same vertex set, say $V$, satisfying certain conditions, there is a perfect matching which creates a Hamilton cycle with each of the matchings in $\cM$.
A related problem, which is motivated by questions in statistical physics, is to compute the maximum over all perfect matchings $M$ on $V$ of the number of cycles that $M$ creates with the given matchings in~$\cM$. 
This has been recently studied when $\cM$ consists of perfect matchings chosen uniformly at random~\cite{GJS:25+}, and we refer to Chapter~2 of the book of Gurau~\cite{guaru:17} for more background.

\subsection{Main results}
Our first main result is a variant of Kotzig's conjecture for random graphs.
The random graph model we consider here is the \defn{binomial random graph $G(n, p)$}, which has $n$ vertices, and where each pair of vertices forms an edge independently with probability $p$. 
Moreover we always assume that $n$ is even.

Since $G(n,p)$ is not regular, it does not have a $1$-factorisation.
However, one can still ask when $G(n,p)$ contains a certain number of perfect matchings, say $k$, such that the union of any two distinct of them induces a Hamilton cycle.
When $k=1$, we are simply requiring the existence of a perfect matching and the threshold is $\log n/n$.
When $k=2$, the problem is equivalent to the existence of a Hamilton cycle (as a Hamilton cycle on an even number of vertices decomposes into two edge-disjoint perfect matchings) and the threshold is again $\log n/n$ by a celebrated result of P\'osa~\cite{posa:76}.
Already for $k=3$, this is not clear anymore.
We show that, for any fixed constant $k$, the result is true above the same threshold.
\begin{theorem}\label{thm:random}
For every fixed $k \in \mathbb{N}$ there exists a constant $C=C(k)$ such that, when $p\ge \frac{C \log n}{n}$, w.h.p.~$G(n,p)$ contains $k$ edge-disjoint perfect matchings $M_1, \dots, M_k$ such that $M_i \cup M_j$ induces a Hamilton cycle for all distinct $i,j \in [k]$.
\end{theorem}

It is tempting to conjecture that the conclusion is already true when $p\ge \frac{\log n+(k-1)\log\log n + \omega(n)}{n}$, or even a hitting time result could hold. We will discuss this further in Section~\ref{sec:remarks}.

Roughly speaking, our approach to prove Theorem~\ref{thm:random} is to find the perfect matchings sequentially by revealing the random edges in $k$ rounds. Assuming that we have already found a collection of perfect matchings such that the union of any two of them induces a Hamilton cycle, we then want to be able to find yet another perfect matching which forms a Hamilton cycle with each of the given ones. Even when the host graph is complete, it does not seem obvious how to achieve this.

Hence the above approach motivates the following natural question: Given a collection of $k$ (edge-disjoint) perfect matchings on a common vertex set, can one find another perfect matching that forms a Hamilton cycle with each of the given matchings?
If $k=2$, then this is easy and in fact for \emph{any} two perfect matchings there exists a third one with the desired property. This can be seen by just adding the edges of the new matching greedily, and in each step, there will be some choice which does not create a forbidden cycle, until the last step.
However, this simple procedure breaks down when $k\ge 3$.
In fact, it is not true anymore that for any given $k$ perfect matchings, there exists another perfect matching forming a Hamilton cycle with each of them, not even when $n$ is arbitrarily large.
For instance, assume $k=3$ and take three matchings $M_1, M_2, M_3$ which are identical on $(n-4)/2$ edges and are pairwise disjoint on the remaining $4$ vertices, which we denote by $v_1,\dots,v_4$ in such a way that $v_1v_2,v_3v_4 \in M_1$.
Suppose by contradiction that there is a perfect matching $M^*$ which creates a Hamilton cycle with each~$M_i$.
Then $M^* \cup M_1$ is a Hamilton cycle and, without loss of generality, we can assume that the cycle traverses the vertices $v_1, \dots, v_4$ in this order.
In particular, there are two non-spanning paths between $v_2$ and $v_3$, and between $v_4$ and $v_1$ which alternate between edges of $M_1$ and edges of $M^*$.
Let $j \in \{2,3\}$ be such that $v_1v_4,v_2v_3 \in M_j$.
Then $M^* \cup M_j$ consists of two cycles, which is a contradiction.

However, our next result shows that, as long as each matching has some edges which belong ``exclusively" to it, we can find the desired additional perfect matching. 

\begin{theorem}\label{thm:abs}
For every fixed $k\in \bN$, there exists a constant $C=C(k)$ such that the following holds. Let $V$ be a set of $n$ vertices and suppose $M_1,\dots,M_k$ are perfect matchings on~$V$, such that
\begin{align}
|M_i\sm (\cup_{j\neq i}M_j)| \ge C \label{exlusiveness}
\end{align}
for each $i\in[k]$.
Then there exists a perfect matching $M^*$ such that $M^*\cup M_i$ forms a Hamilton cycle for each $i\in [k]$.
\end{theorem}

We will prove Theorem~\ref{thm:abs} using the absorbing method which allows to ``plan ahead'' and prepare for a situation where the greedy algorithm gets stuck.
As already discussed above, condition~\eqref{exlusiveness} is not needed when $k=2$, but cannot be omitted when $k \ge 3$.

\cref{thm:abs} establishes the existence of one perfect matching $M^*$ which forms a Hamilton cycle with each of the given matchings.
The next natural question is to count the number of possibilities for $M^*$. To motivate our result, recall that the number of perfect matchings in $K_n$ is $\frac{n!}{2^{n/2}(n/2)!}=\left(\sqrt{2} \pm o(1)\right) \cdot \left((n/e)^{n/2}\right)$. 
Suppose now that a perfect matching $M_1$ is given and we want to count the number of perfect matchings $M_2$ of $K_n$ such that $M_1\cup M_2$ forms a Hamilton cycle.
The answer is easily seen to be $\frac{2^{n/2} \cdot (n/2)!}{2 \cdot (n/2)}= \left(\sqrt{\pi} \pm o(1)\right) \cdot \left(n^{-1/2}(n/e)^{n/2}\right)$.
In other words, if $M_2$ was chosen uniformly at random from all perfect matchings, then the probability that it forms a Hamilton cycle together with $M_1$ is $\left(\sqrt{\frac{\pi}{2}} \pm o(1)\right) \cdot n^{-1/2}$.
Heuristically, one might then expect that given $k$ ``independent'' (say, edge-disjoint) matchings, the probability that a uniformly random perfect matching forms a Hamilton cycle with each of them is $\left(\left(\frac{\pi}{2}\right)^{k/2} \pm o(1)\right) \cdot n^{-k/2}$.
When $k=2$, that is, in the case where the greedy algorithm never gets stuck, Kim and Wormald~\cite{KW:01} proved that this is indeed correct. By applying this iteratively, they also obtained the asymptotic number of matchings $M_1, \dots, M_k$ such that $M_i \cup M_j$ forms a Hamilton cycle for each $ij \in E(G)$ with $G$ being a fixed graph on $[k]$ of degeneracy~$2$. 
However, when $k>2$, their method does not work anymore, at least not without significant new ideas that allow for some ``planning ahead''.
We confirm the heuristic argument for any constant $k$, up to the multiplicative constant factor.

\begin{theorem}\label{thm:count}
Let $k\in \bN$ be fixed and $n$ sufficiently large. Let $V$ be a set of $n$ vertices and suppose $M_1,\dots,M_k$ are edge-disjoint perfect matchings on~$V$.
Then the number of perfect matchings $M^*$ such that $M^*\cup M_i$ forms a Hamilton cycle for each $i\in [k]$, is $\Theta_k\left(n^{-k/2}(n/e)^{n/2}\right)$.
\end{theorem}

The condition that the given matchings $M_1,\dots,M_k$ are edge-disjoint cannot be dropped. For instance, if they would be identical, then the number of choices for $M^*$ is $\Theta\left(n^{-1/2}(n/e)^{n/2}\right)$ as discussed above.
In fact, we will prove a more general version where we can relax the edge-disjointness condition to a mild overlap (cf.~\cref{thm:count robust}).
For its proof, we use a variety of methods, including a random process analysis, the absorption method (via \Cref{thm:abs}), a variant of the entropy method and the switching method. 
We then use this stronger result to deduce Theorem~\ref{thm:random} via the recent breakthroughs on the expectation threshold conjecture (cf.~\cref{sec:proof_thm_random}), where our goal is to show that the auxiliary hypergraph whose hyperedges correspond to the possible matchings $M^*$ is suitably spread.

Let us briefly discuss the methods we use. As already mentioned, the proof of the existence result employs the absorbing method, which is a powerful method that in general is applied to turn approximate solutions into exact solutions. While the greedy algorithm can potentially get stuck, this will only happen towards the end of the process. Very roughly speaking, we will overcome this by planting some special edges in the beginning that will belong to the final matching and will yield some flexibility to resolve any final situation where the greedy algorithm gets stuck. We hope that our ideas will shed some light on potential ways of attack on Kotzig's conjecture using such kind of methods.

For the counting result, a natural approach is to determine the number of choices the greedy algorithm has in each step, and multiply them together. Unfortunately, this sequence is not deterministic, but depends on the history of already chosen edges. A commonly used arsenal of methods for such counting problems of combinatorial configurations is the following: For lower bounds, to use a randomized greedy algorithm whose behaviour can be analysed through martingale inequalities. Here, the goal is to show that, typically, the number of choices in each step follows some deterministic trajectory, up to small error terms, and that, by stopping the process early enough, a typical outcome of the produced approximate solution can be completed (using the absorbing method, say). For upper bounds, the entropy method has proven itself as a very versatile tool. Here, the main idea is that the number of solutions is directly linked to the entropy of a uniformly random solution. One can then apply known facts about entropy like the chain rule to analyse this in a step-by-step manner akin to the random greedy algorithm.
While our approach is certainly inspired by these methods, the main difficulties in our case are quite different from previous applications. Perhaps most notably, the key random variable which controls the number of choices in each step, which is the so-called total overlap (cf. Definition~\ref{def:overlap}), is not concentrated. In fact, one heuristically expects it to be bounded by a constant. While this means one should be able to show that it typically never exceeds $\Theta(\log n)$, this would not be enough to obtain a counting result that is tight up to constant factors.

Hence, one of the key features of our counting proof is that we analyse the total contribution of this variable over the whole process. Fortunately, it turns out that it is enough to control the expectations. For this, we exploit a ``self-correcting'' behaviour. Very roughly speaking, this means that if the total overlap exceeds some constant threshold, then there will be a negative drift that pushes it down again. This phenomenon is also crucial for the strengthening of our counting result where the matchings can have a mild overlap at the start (which in turn is essential for the result on random graphs). For this, we need to show that the total overlap decreases quickly and then stays below the constant threshold. 

One caveat with using the entropy method is that we have to work with a solution chosen uniformly at random, which is a difficult distribution to handle due to the lack of independence. In many applications, one can bypass this issue with the following conditioning trick: fix any solution (say a perfect matching with the desired property) and reveal its elements (the edges) in a random order. If one can efficiently control the expected number of choices for each step (with respect to the previously appearing elements of the fixed solution), then this may already yield the desired bound.
However, in our case, this does not seem tractable, due to the rather complicated and ``long-range'' possibilities that can make an edge unavailable. (This is slightly reminiscent of the counting problem for so-called ``high-girth'' designs, where so far the conjectured upper bounds could not be established.) Thus, since we really have to exploit the properties of a uniformly random solution, a final key ingredient in our analysis is a statistical result which says that every available edge is roughly equally likely to be contained in a uniformly random solution. This uses a cascade of several switching operations.

\section{Notation, organisation and basic concepts}

\subsection{Notation}
We let $[n]$ denote the set $\{1, \dots, n\}$ and write $\log(x)$ to mean $\log_e(x)$.

We use standard graph theory terminology.
Given a matching $M$ we denote by $V(M)$ the set of vertices covered by~$M$.
We denote an ordered matching with $m$ edges by $(e_1,\dots,e_m)$.

We say that an event holds \defn{with high probability} (w.h.p.) if the probability that it holds tends
to~$1$ as the number of vertices $n$ tends to infinity.

For $a, b, c \in (0, 1]$, we write $a \ll b \ll c$ in our statements to mean that there are increasing
functions $f, g : (0, 1] \to (0, 1]$ such that whenever $a \le f(b)$ and $b \le g(c)$, then the subsequent
result holds. Moreover, when using the Landau symbols $O(\cdot), \Omega(\cdot), \Theta(\cdot)$, subscripts denote variables that the implicit constants may depend on.
Similarly, $o_D(\cdot)$ stands for a function which tends to $0$ as $D \to \infty$.

\subsection{Organisation of the paper}
We introduce some basic concepts in the next subsection.
We prove \cref{thm:abs} in \cref{sec:existence} and (the strengthening of) \cref{thm:count} in \cref{sec:counting}.
The proof of \cref{thm:count} will require an additional result which we prove in \cref{sec:switiching}.
We prove \cref{thm:random} in \cref{sec:proof_thm_random}, where we also discuss the relevant concepts around the expectation threshold conjecture.
We then discuss analogous results in the bipartite setting in \cref{sec:bipartite} and finish with some concluding remarks in \cref{sec:remarks}.

\subsection{Basic concepts}
\label{subsec:basic_concepts}
We introduce some basic concepts which will be used throughout the paper.
\begin{defin}[$k$-configuration]
A \defn{$k$-configuration} is a set of $k$ (not necessarily edge-disjoint) perfect matchings on a common vertex set.
The \defn{order} of a $k$-configuration is the size of its vertex set.
We say that a pair of vertices forms an \defn{available} edge if the edge is not contained in any of the given perfect matchings. 
\end{defin}

We will often associate colours with the $k$ matchings. For two generic matchings, we will use blue and green.
Given a $k$-configuration, we aim to find a matching $M^*$ that forms a Hamilton cycle with each of the given matchings. We will always refer to this matching as the \defn{red} matching. 

\begin{figure}[htp]
\begin{tabular*}{\textwidth}{@{\extracolsep{\fill}}cccc}  
  \includegraphics[width=.22\textwidth]{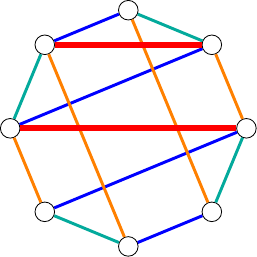}
& \includegraphics[width=.22\textwidth]{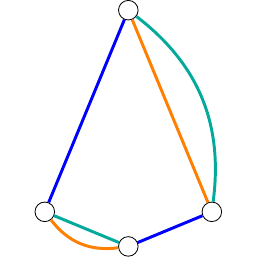} 
& \includegraphics[width=.22\textwidth]{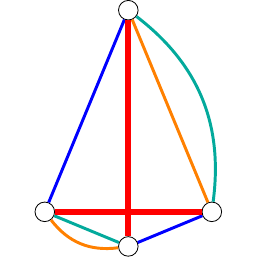}  
& \includegraphics[width=.22\textwidth]{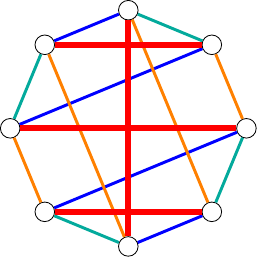} \\
  (A) & (B) & (C) & (D)\\
\end{tabular*}
    \caption{(A): A $3$-configuration $\cM$ of order $8$ and a partial red matching $M$ of size $2$. (B): The reduced configuration of $\cM$ with respect to $M$ (note that each colour still induces a perfect matching and that we have created double coloured edges). (C): A red perfect matching $M'$ of the reduced configuration.
    (D): $M \cup M'$ is a red perfect matching of $\cM$.}
    \label{fig:reduced_configuration}
\end{figure}

We will construct the red matching in an iterative way.
Suppose we are given a \defn{partial red} matching. When speaking of a partial red matching, we always assume that it is properly partial, that is, not a perfect matching yet, and we also assume that the union of this red matching with any other colour is still cycle-free. (Once such a cycle would be closed, it would be impossible to extend this red matching to a full red matching.)
We would like to add a red edge, while keeping the property that the union of the new red matching with any other colour is still cycle-free. We encode the edges which we are not allowed to pick as the new red edge through the following auxiliary graph (cf.~(B) in \cref{fig:reduced_configuration}).

\begin{defin}[Reduced configuration]
\label{def:reduced_conf}
Given a $k$-configuration $\cM$ and a partial red matching~$M$, we define the \defn{reduced configuration} as follows. Its vertex set is the set of vertices not covered by an edge of $M$, which we denote by $U$. Then for each colour, say blue, and each vertex $u\in U$, we consider the path starting at $u$ with edges alternating blue and red and we let $v$ be the other endpoint. (Clearly $v \neq u$ and $v \in U$.) We then insert $uv$ as a blue edge in the reduced configuration (cf.~\cref{fig:reduced_configuration}).
\end{defin}

Observe that the reduced configuration is still a $k$-configuration, in the sense that each colour induces a perfect matching, and that an edge of the reduced configuration may be assigned multiple colours. Moreover, if the $k$-configuration has an edge $e$ of some colour, say blue, and none of its endpoints is covered by a red edge, then $e$ is still an edge of colour blue in the reduced configuration.
The reduced configuration conveniently allows us to ``forget'' the red edges that we have already inserted: Indeed, the collection of red edges inserted so far can be completed to a red perfect matching if and only if its reduced configuration admits a red perfect matching (cf.~\cref{fig:reduced_configuration}).

The example in \Cref{fig:reduced_configuration} shows that, even if we start from a configuration with pairwise edge-disjoint perfect matchings, we may create edges with multiple colours in the reduced configuration.
However, for our arguments, it is essential to control how the number of edges which get more than one colour evolves in the corresponding reduced configurations after adding a new red edge, which we capture through the following definition.

\begin{defin}[Overlap]
\label{def:overlap}
    For a $k$-configuration $\cM:=\{M_1,\dots,M_k\}$ and two colours $i, j$, we define the \defn{overlap} of colours $i$ and $j$ as $R_{ij}(\cM):=|M_i \cap M_j|$. 
    Moreover, we define the \defn{total overlap} as $R(\cM):=\sum_{1 \le i < j \le k} R_{ij}(\cM)$.
\end{defin}

In order to understand the behaviour of the total overlap, it is enough to consider a single pair of colours at the time.
Therefore, fix any two colours, say $i$=blue and $j$=green, and note that the union of the blue and green matchings decomposes into alternating blue-green cycles and double coloured edges.
Then observe that there are five different ways we can add a red edge, depending on the components (of the blue-green graph) its endpoints belong to:
\begin{enumerate}[label=\rm{(\alph*)}]
\item \label{case_A} both endpoints belong to the same cycle $C$ and their distance is odd;
\item \label{case_B} both endpoints belong to the same cycle $C$ and their distance is even;
\item \label{case_C} its endpoints belong to different cycles $C_1$ and $C_2$;
\item \label{case_D} one endpoint belongs to a cycle $C$ and the other one to a double edge $f$;
\item \label{case_E} both endpoints belong to double edges, say $f_1$ and $f_2$.
\end{enumerate}
For each case, we can easily describe the changes to the blue-green cycle structure of the reduced configuration.
In case~\ref{case_A}, $C$ splits into two shorter cycles; in particular, if the endpoints of the red edge are at distance exactly three, then at least one of these cycles would be a (new) double edge and $R_{ij}(\cdot)$ would increase by $1$ or $2$. 
In case~\ref{case_B}, $C$ becomes a cycle whose length is by two shorter than that of $C$; in particular, if the length of $C$ was $4$, the shorter cycle would be a (new) double edge and $R_{ij}(\cdot)$ would increase by $1$.
In case~\ref{case_C}, the two cycles $C_1$ and $C_2$ merge into a new longer cycle, and $R_{ij}(\cdot)$ does not change.
In case~\ref{case_D}, $C$ and $f$ merge into a cycle; in particular, $R_{ij}(\cdot)$ decreases by $1$.
In case~\ref{case_E}, $f_1$ and $f_2$ merge into a double edge; in particular, $R_{ij}(\cdot)$ decreases by $1$.
In particular, we note that the addition of a new red edge changes $R_{ij}(\cdot)$ by at most $2$, and thus the total overlap changes by at most $2 \binom{k}{2} \le k^2$.
In order to keep track of $R_{ij}(\cdot)$ over time, we make the following definition.

\begin{defin}[Good and bad edges]
\label{def:good_edges}
    Let $\cM$ be a $k$-configuration and $e \not \in \cM$.
    Let $\cM'$ be the reduced configuration of $\{e\}$.
    Given two colours, say $i,j$, we say that $e$ is \defn{good} (resp.~\defn{bad}) for the colours $i$ and $j$ if $R_{ij}(\cM') < R_{ij}(\cM)$ (resp.~$R_{ij}(\cM') > R_{ij}(\cM)$).
    We remark that the available edges whose addition does not change $R_{ij}(\cdot)$ are neither good nor bad.
\end{defin}

The discussion above immediately implies the following proposition.

\begin{prop}
\label{prop:counting_double_edges}
    Let $\cM$ be a $k$-configuration of order $n$ and $i, j$ be two colours.
    Set $R_{ij}:=R_{ij}(\cM)$.
    Then the following properties hold.
    \begin{enumerate}[label=\rm{(P\arabic*)}]
        \item \label{prop_P_1} The number of good edges is at least $(R_{ij}/2-k/2)n$ and the number of bad edges is at most $n$.
        \item \label{prop_P_2} Let $e \not\in \cM$ and $\cM'$ be its reduced configuration. If $e$ is a good edge then $R_{ij}(\cdot)$ decreases by $1$ and if $e$ is a bad edge then $R_{ij}(\cdot)$ increases by $1$ or $2$.
    \end{enumerate}
\end{prop}
\proof
    Let $S$ be the set of vertices covered by the double edges, i.e.~those in $M_i \cap M_j$.

    The good edges are those corresponding to cases~\ref{case_D} and~\ref{case_E}, i.e.~the available edges which are incident to a vertex of $S$.
    Their number is at least $|S|(n-|S|)+\binom{|S|}{2} - k (n/2) \ge (|S|/4-k/2)n= (R_{ij}/2-k/2)n$, where we used $|S| \le n$.
    
    Similarly, the bad edges are those corresponding to case~\ref{case_A} (with the endpoints of the red edge being at distance three on the cycle) and case~\ref{case_B} (with the cycle having length four).
    Therefore, for each vertex, there exist at most two bad red edges incident to this vertex, and thus the number of bad edges is at most $n$.
    This proves Item~\ref{prop_P_1}.
    
    Item~\ref{prop_P_2} follows directly from the discussion above of cases~\ref{case_A},~\ref{case_B},~\ref{case_D} and~\ref{case_E}.
\endproof

Finally, we prove an easy observation.

\begin{observation}
\label{obs:random_stays_random}
    Let $\cM$ be a $k$-configuration of order $n$, fix $m \le n/2$ and let $M^*=(e_1^*,\dots,e_m^*)$ be an ordered red matching of size $m$ chosen uniformly at random.
    Then, for each $t \le m$, conditioned on $e_1^*,\dots,e_{t-1}^*$, the ordered matching $(e_t^*,\dots,e_m^*)$ is distributed as a uniformly chosen ordered red matching of size $m-t+1$ of the reduced configuration of $\{e_1^*,\dots,e_{t-1}^*\}$.
\end{observation}

\proof
    Let $\cX$ be the set of ordered red matchings of size $m$ of $\cM$.
    Let $(f_1,\dots,f_m) \in \cX$ and condition on $e_i^*=f_i$ for each $i \le t-1$.
    Denote by $\cM'$ the reduced configuration of $\{f_1,\dots,f_{t-1}\}$.
    Then
    \[
        \mathbb{P}\Bigg[\bigcap_{i=t}^m \{e_i^*=f_i\} \bigg| \bigcap_{i=1}^{t-1} \{e_i^*=f_i\}\Bigg] = \frac{\prob{\cap_{i=1}^m \{e_i^*=f_i\}}}{\prob{\cap_{i=1}^{t-1} \{e_i^*=f_i\}}} = \frac{1/|\cX|}{|\cA|/|\cX|} = \frac{1}{|\cA|} \, ,
    \]
    where $\cA$ denotes the set of matchings of $\cX$ with $e_i^*=f_i$ for each $i \le t-1$.
    We conclude by observing that $\cA$ has the same cardinality as the set of the red matchings of size $m-t+1$ of $\cM'$.
\endproof

\section{Existence -- Proof of Theorem~\ref{thm:abs}}
\label{sec:existence}
We first discuss the proof of \cref{thm:abs}.
Suppose that we are given a $k$-configuration and we would like to construct a red perfect matching.
We start by observing that we can greedily build a partial red matching of size at least $\frac{n-k-1}{2}$.
Indeed, suppose that we have already chosen some red edges, consider the reduced configuration and let $n'$ be its number of vertices.
The union of the reduced edges has size at most $\frac{n'}{2}\cdot k$ and, as long as $n' 
> k+1$, we have $\frac{n'}{2}\cdot k<\binom{n'}{2}$ and thus there is a choice for the next red edge.
The main part of the proof is to show that the reduced configuration obtained when there are only few uncovered vertices ($k$ or $k+1$, depending on the parity of $k$) can be solved, i.e.~that it has a red perfect matching. Note that not all configurations are solvable: For example, if we take three perfect matchings on four vertices whose union forms a $K_4$, there is no red perfect matching. 
Therefore we employ the absorbing method, building a flexible structure in the beginning which will help us to deal with any potential leftover.

To motivate our approach, observe that if the $k$ matchings of a $k$-configuration are identical, then a red perfect matching clearly exists.
Hence, if we could choose some red edges so that in the reduced configuration the $k$ colours become identical, we would be done.
When $k=2$, this can be easily achieved.

\begin{prop}\label{prop:two colours}
Given a $2$-configuration, there exists a partial red matching such that in the reduced configuration the two matchings are identical.
\end{prop}

\proof
Let the two colours be blue and green.
Recall that the union of any two perfect matchings forms a spanning subgraph whose components are alternating blue-green (even) cycles or double edges. For the latter case, we do not have to do anything, so suppose we have an even cycle with edges alternating blue and green. If the length is at least $6$, then taking a red edge between two vertices at distance $3$ on the cycle splits off a double-coloured edge and a shorter cycle in the reduced configuration. Finally, if the cycle has length $4$, then choosing one of the diagonals as a red edge leaves a single double-coloured edge in the reduced configuration. By repeatedly choosing red edges in this way, we find a partial red matching with the desired property.
\endproof

Once two colours are identical, they will remain identical, regardless of which additional red edges we add. Hence, a natural approach is to equalize the colours step by step. Note that this is not as easy as for $k=2$: When $k>2$, we might not be able to add the red edges as in the proof of Proposition~\ref{prop:two colours} since these might be edges of a third colour.

\begin{defin}[Equalizing matching]
Given a $k$-configuration $\cM$ and two colours, say blue and green, a partial red matching is called \defn{blue-green equalizing} if in the reduced configuration the blue and green matchings are identical. It is called \defn{totally equalizing} if in the reduced configuration all the $k$ matchings are identical.
\end{defin}

Observe that Proposition~\ref{prop:two colours} can be rephrased equivalently by saying that any $2$-configuration has a totally equalizing partial red matching.
We remark once more that we only seek a partial red matching which is properly partial, meaning that its reduced configuration is still non-empty. This is because we will find such equalizing matchings in different subconfigurations of a given configuration, and in the end we want to find a solution which ``connects'' these different parts.

We now define our absorber.

\begin{defin}[Absorber]
    Let $\cL$ be a $k$-configuration, and blue and green be two colours.
    Let $\cA$ be a $k$-configuration, which is vertex-disjoint from $\cL$ (but has the same colours).
    Then $\cA$ is called a \defn{blue-green equalizer} for $\cL$ if both $\cA$ and $\cA\cup \cL$ have a blue-green equalizing red matching.
    Moreover $\cA$ is called an \defn{absorber} for $\cL$ if both $\cA$ and $\cA\cup \cL$ have a totally equalizing red matching.
\end{defin}

We remark that, while $\cA$ and $\cL$ are vertex-disjoint, the equalizing matching of $\cA \cup \cL$ is still allowed to contain edges between $\cA$ and $\cL$.
\medskip

\subsection{Existence of absorbers}
Given a $k$-configuration, we will construct an absorber by iteratively building equalizers for two colours at the time.
Therefore we start by showing how we can equalize two colours.

\begin{lemma}\label{lem:equalizer}
Given a $k$-configuration $\cL$ of order $m$ and two colours, say blue and green, there exists a blue-green equalizer of order $2m$.
\end{lemma}

\begin{figure}[t!]
\centering
    \includegraphics[scale=.6]{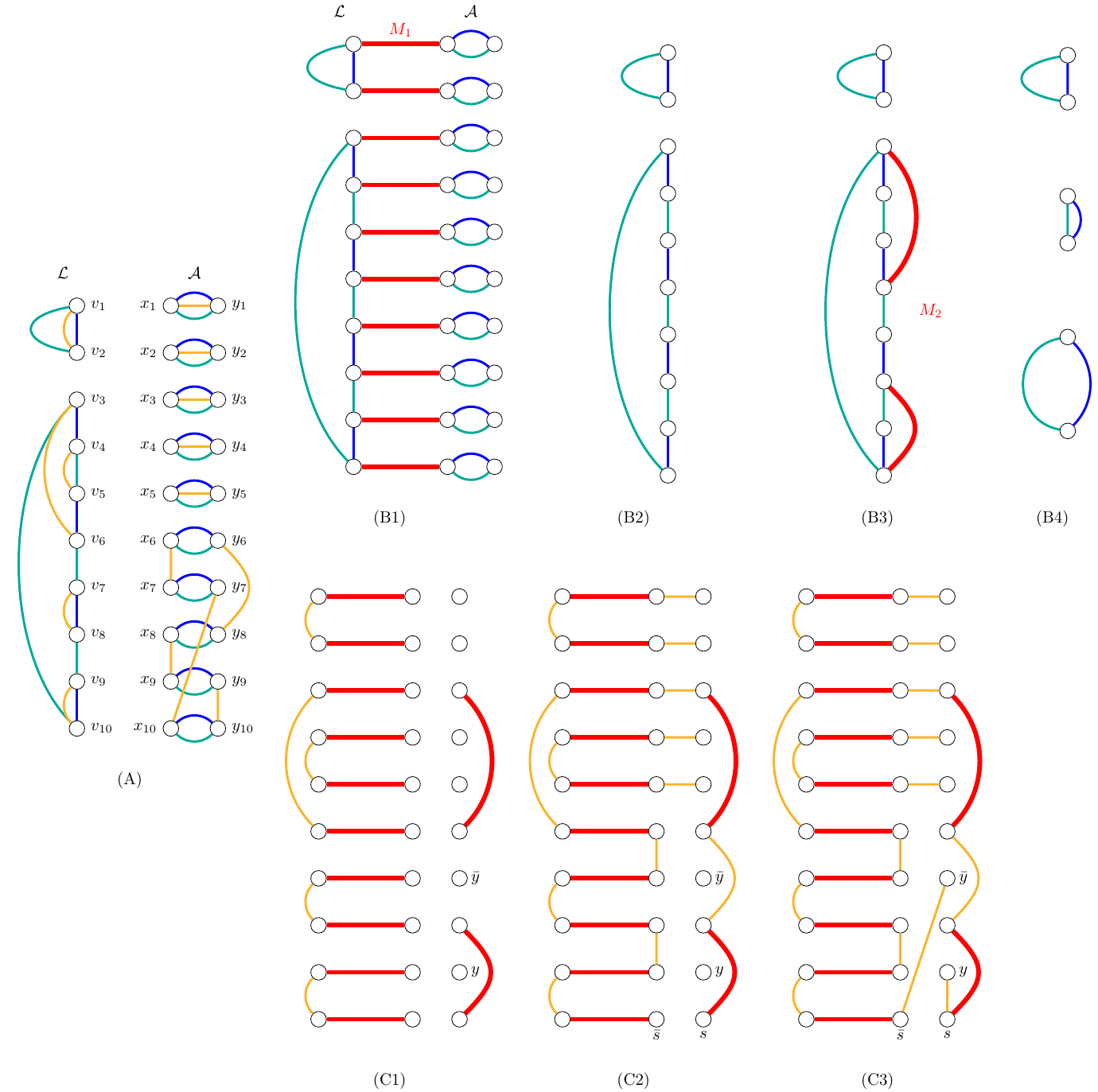}
    \caption{(A): A $3$-configuration $\cL$ on $10$ vertices $\{v_i\}$ and its blue-green equalizer $\cA$ with vertices $\{x_i\} \cup \{y_i\}$: The empty matching is blue-green equalizing for $\cA$ and $M_1 \cup M_2$ is blue-green equalizing for $\cA \cup \cL$.\\
    (B): We show how we add the blue and green matching on $Z$ and how we choose the red edges in $M_1$ and $M_2$ in the proof of \cref{lem:equalizer}. (B1): We first add double blue-green coloured edges $x_iy_i$ and set $M_1:=\{v_ix_i\}$. (B2): The blue and green edges of the reduced configuration of $M_1$. (B3): We find a red matching $M_2$ which equalizes blue and green. (B4): Indeed the blue and green edges of the reduced configuration of $M_1 \cup M_2$ are identical.\\
    (C): We show how we add the yellow matching on $Z$ in the proof of \cref{lem:equalizer}. (C1): Once $M_1$ and $M_2$ have been fixed, we arbitrarily choose two distinct vertices $y, \bar{y}$ not covered by a red edge. (C2): We can greedily add a yellow matching covering all but $4$ vertices such that $y, \bar{y}$ are not covered and we do not close a red-yellow cycle. (C3): We add the yellow edges $sy, \bar{s}\bar{y}$ which completes a yellow perfect matching.}
    \label{fig:equalizing}
\end{figure}

\proof
Let $v_1,\dots,v_m$ be the vertices of $\cL$ and set $Z:=\{x_1, \dots, x_m, y_1, \dots, y_m\}$ to be a set of $2m$ new vertices. We will construct a blue-green equalizer $\cA$ with $V(\cA)=Z$ (\cref{fig:equalizing} illustrates the various steps).

First we take the edges $x_iy_i$ as blue-green double coloured edge for all $i \in [m]$, and denote by $E_1$ their union. Note that these edges give a blue and a green perfect matching in $V(\cA)$ and thus we will not add any more edges of these two colours. Moreover, observe that, no matter which edges we add for the other colours, $\cA$ has a blue-green equalizing red matching, namely the empty matching. 

Now, for each colour different from blue and green, we have to add a perfect matching of $Z$ in this colour in order to turn $\cA$ into a $k$-configuration.
Moreover, we want to choose them in such a way that we can then find a blue-green equalizing red matching $M$ in $\cA \cup \cL$.
In order to take advantage of \cref{prop:two colours}, we will achieve this in the opposite order: We first define which edges we want for $M$ and then show that it is still possible to add the other colours in a ``compatible'' way, meaning that the union of $M$ and any of these matchings is cycle-free (where we also forbid double edges).

We first take $v_ix_i$ as red edge for all $i\in [m]$, and denote by $M_1$ their union (cf.~drawing~(B1) in \cref{fig:equalizing}). Consider the reduced configuration of $M_1$ with respect to colours blue and green (cf.~drawing~(B2) in \cref{fig:equalizing}).
Then its vertex set is $Y=\Set{y_1,\dots,y_m}$ and, by Proposition~\ref{prop:two colours}, there exists a partial blue-green equalizing red matching $M_2$ on $Y$ (cf.~drawings~(B3) and~(B4) in \cref{fig:equalizing}).
In particular, $M_1 \cup M_2$ is a blue-green equalizing red matching for $\cA \cup \cL$.

Now we move to the other colours.
Since the edges of $\cA$ are allowed to receive more than one colour, we can simply focus on one such colour, say yellow.\footnote{In some sense, in this part of the proof yellow and red swap their roles. There are some red edges already fixed and we have to choose a yellow perfect matching while not creating alternating cycles.}
Since $M_2$ is a properly partial matching of $Y$, there are at least two vertices of $Y$, say $y$ and $\bar{y}$, which are not covered by a red edge (cf.~drawing~(C1) in \cref{fig:equalizing}).
We start by covering all but $4$ vertices of $Z$ with a partial yellow matching, while not using the vertices $y$ and $\bar{y}$.
This can be done greedily (cf.~drawing~(C2) in \cref{fig:equalizing}).
Indeed, let $\tilde{Z}$ be the set of vertices of $Z$ which are not (yet) covered by a yellow edge, and suppose $y, \bar{y} \in \tilde{Z}$ and $|\tilde{Z}| >4$.
Then, take any $z \in \tilde{Z} \setminus\{y,\bar{y}\}$, consider (if it exists) the red-yellow alternating path starting from $z$ and let $z'$ be its other endpoint.
Then there exists $w \in \tilde{Z} \setminus\{y,\bar{y},z\}$ with $w \neq z'$ and we can add $zw$ as a yellow edge.
(If the alternating path does not exist, we can take as $w$ any vertex in $\tilde{Z} \setminus\{y,\bar{y},z\}$.)
At the end of the greedy process, there are exactly $4$ vertices which are not covered by a yellow edge: these are $y, \bar{y}$ and two more, say $s$ and $\bar{s}$.
We can finish by adding $ys$ and $\bar{y}\bar{s}$ as yellow edges as, since neither $y$ nor $\bar{y}$ is covered by a red edge, this addition cannot close a red-yellow cycle 
(cf.~drawing~(C3) in \cref{fig:equalizing}).
We let $E_2$ be the union of all the edges added in this step (among all colours).

Observe that the configuration on $Z$ with edges $E_1 \cup E_2$ gives the desired equalizer (cf.~drawing~(A) in \cref{fig:equalizing}): Indeed, $E_1 \cup E_2$ is the union of $k$ perfect matchings of $Z$, the empty matching is blue-green equalizing for $\cA$, and $M_1 \cup M_2$ is blue-green equalizing for $\cA \cup \cL$.
\endproof

We can then show the existence of absorbers.

\begin{lemma}\label{lem:abs}
Any $k$-configuration $\cL$ of order $m$ has an absorber.
\end{lemma}

\proof
We prove the statement by induction on $k$.
If $k=2$, the statement is true by \cref{lem:equalizer}.
Suppose it is true for $k-1$ and let $\{c_1,\dots,c_k\}$ be the set of the colours of the matchings in $\cL$.
By \cref{lem:equalizer}, there exists an equalizer $\cA$ for $\cL$ for the colours $c_1$ and $c_2$.
Let $M_1$ (resp.~$M_2$) be the $c_1$-$c_2$ equalizing red matching of $\cL \cup \cA$ (resp.~$\cA$).

Let $\cL_1$ be the reduced configuration of $\cL \cup \cA$ with respect to the matching $M_1$. 
Then, in $\cL_1$, the matchings of colour $c_1$ and $c_2$ are equal and we can consider them as being the same colour. 
In particular, by induction, there exists an absorber $\cA_1$ of $\cL_1$ (for the colours $c_2, \dots, c_k$ and thus $c_1$ as well).
We let $M_1^1$ (resp.~$M_2^1$) be the totally equalizing red matching of $\cL_1 \cup \cA_1$ (resp.~$\cA_1)$.
Similarly, let $\cL_2$ be the reduced configuration of $\cA$ with respect to the matching $M_2$. 
Then, by induction, there exists an absorber $\cA_2$ of $\cL_2$ (for the colours $c_2, \dots, c_k$ and thus $c_1$ as well).
We let $M_1^2$ (resp.~$M_2^2$) be the totally equalizing red matching of $\cL_2 \cup \cA_2$ (resp.~$\cA_2)$.

Let $\cA':=\cA \cup \cA_1 \cup \cA_2$ (notice the union is disjoint).
Then $\cA'$ is a $k$-configuration and is an absorber for $\cL$.
Indeed, $M_1 \cup M_1^1 \cup M_2^2$ is a totally equalizing red matching for $\cL \cup \cA'$ while $M_2 \cup M_2^1 \cup M_1^2$ is a totally equalizing red matching for $\cA'$.
\endproof
\medskip

\subsection{Finding absorbers inside configurations.}
While we have shown that absorbers exist, we still have to justify why we can indeed ``find'' them inside a given initial $k$-configuration. By ``finding'', we mean that we can choose a few red edges such that the desired absorber appears inside the reduced configuration. 
This is the step where we crucially need that every matching has some exclusive edges.

We first show that taking at most two red edges suffices to force a fixed edge to appear in the reduced configuration and in the desired colour, say blue, provided there are enough edges which are exclusively blue.

\begin{lemma}\label{lem:embed}
Let $\cM$ be a $k$-configuration on $V$ and let $Z \subseteq V$. Suppose that there are at least $(|Z|+2)\cdot k+1$ edges with both endpoints in $V \setminus Z$ which are exclusively blue.
Let $x,y\in Z$ such that the blue neighbours of both $x$ and $y$ belong to $V \setminus Z$. 
Then there exists a partial red matching $M$ of size $2$ such that its edges lie outside $Z$ and its reduced configuration has the following property:
The only new edge with both endpoints in $Z$ is $xy$ and its colour is blue, and all the others stay the same.
\end{lemma}

\proof
    Let $x'$ and $y'$ be the other endpoints of the blue edges adjacent to $x$ and $y$, respectively.
    By assumption $x',y' \not \in Z$.
    Set $T:=N\left(Z \cup \{x',y'\}\right)$, i.e.~$T$ is the set of vertices $t$ such that $tw$ is an edge of $\cM$ (of any colour) for some $w \in Z \cup \{x',y'\}$.
    Since $|T| \le (|Z|+2)k$, there exist $x'',y'' \in V \setminus (Z \cup T)$ such that the edge $x''y''$ is exclusively blue.

    We set $e:=x'x''$ and $f:=y'y''$, and claim that $M:=\{e,f\}$ is a feasible red matching.
    Indeed, clearly it does not close any alternating red-blue cycle.
    For any other colour, say green, since $e,f$ are available edges in $\cM$, none of them is green and thus adding $M$ does not create any double red-green edge.
    Suppose now there is an alternating red-green cycle using both these two edges.
    Then, since there are no other red edges, the cycle has length $4$ and it is either $x'x''y''y'$ or $x'y''y'x''$.
    The first cannot happen as otherwise $x''y''$ would be a green edge, while it is exclusively blue by assumption.
    The latter cannot happen either as otherwise $x'y''$ would be a green edge, contradicting the definition of $T$.

    Observe that in the reduced configuration of $M$, the edge $xy$ is blue.
    We now show that no other new edge is created in $Z$ and that all the edges with both endpoints in $Z$ stay the same.
    Suppose $vw$ is a new edge of the reduced configuration of $M$ with $v,w \in Z$, and let $c$ be its colour. 
    Then there is a $c$-red alternating path with endpoints $v$ and $w$, using some of the added red edges $e$ and $f$.
    By definition of $T$ and the choice of $x''$ and $y''$, among the endpoints of $e$ (resp.~$f$) only $x'$ (resp.~only $y'$) can be adjacent in colour $c$ to $v$ or $w$, and if this is the case then it cannot be adjacent to both $v$ and $w$.
    Therefore, without loss of generality, we can assume that the $c$-red alternating path is $vx'x''y''y'w$.
    Since $x''y''$ is exclusively blue, then the colour $c$ is blue and the edge $vw$ is $xy$, as wanted.
    Finally, since $M$ lies outside $Z$, all the edges with both endpoints in $Z$ stay the same.
\endproof

Next, we apply the previous lemma several times in order to find the desired absorber.

\begin{lemma}\label{lem:find}
For all integers $k,m \ge 2$, there exists a constant $C=C(k,m)$ such that the following holds for any $k$-configuration $\cA$ of order $m$.
Suppose that $\cM$ is a $k$-configuration with the property that, for each of the $k$ colours, $\cM$ has at least $C$ edges of that colour only.
Then there exists a red matching $M^*$ such that the reduced configuration contains~$\cA$.
\end{lemma}

\proof
    Set $C:=k \cdot [(k+3)m+3]$.
    Let $\cE:=\{e_1,\dots,e_{\ell}\}$ be the multiset of edges of $\cA$ and observe that $\ell = (km)/2$.
    Given the edge $e_i$, denote its endpoints by $a_i$ and $b_i$, and its colour by $c_i$. 
    Let $V$ be the vertex set of $\cM$, and fix an independent subset $Z \subseteq V$ of size $m$ and any bijection $\phi:V(\cA) \to Z$. 
    (Such a set $Z$ exists as $\cM$ is the union of $k$ perfect matchings and has at least $2C$ vertices.)
    We show that for each $i \in [\ell]$
    \begin{enumerate}[label = (Q)]        
        \item \label{eq:Q} there exists a red matching $M^*_i$ of size $2i$ such that its edges lie outside $Z$ and the edges of its reduced configuration with both endpoints in $Z$ are precisely $\phi(a_j)\phi(b_j)$ in colour $c_j$ for each $j \in [i]$.
    \end{enumerate}
    Then we can set $M^*:=M^*_{\ell}$ and the lemma follows.

    The case $i=1$ of~\ref{eq:Q} follows from Lemma~\ref{lem:embed} on input $Z$, $x=\phi(a_1)$, $y=\phi(b_1)$ and colour~$c_1$.
    (Observe that since $Z$ is an independent set, the neighbours of $x$ and $y$ in colour $c_1$ are outside~$Z$.)
    Now suppose that~\ref{eq:Q} holds for some $i<\ell$ and let $M^*_i$ be the red matching.
    We want to apply \cref{lem:embed} to the reduced configuration of $M^*_i$ (which still contains the set $Z$).
    For that, for each colour, say blue, we need to count the number of its edges which are exclusively blue.
    So fix an edge of $\cM$ which lies outside $Z$ and was exclusively blue at the beginning, say $f:=vw$.
    The edge remains exclusively blue in the reduced configuration of $M^*_i$ unless one of the following circumstances happens.
    First, it could be that one of the endpoints of $f$ belongs to $M^*_i$.
    Since each vertex in $M^*_i$ is adjacent to at most one (exclusively) blue edge, this can happen for at most $|V(M^*_i)|$ such edges.
    Second, it could be that $f$ gets an additional colour, i.e.~there are another colour, say green, and an alternating green-red path with endpoints $v$ and $w$.
    The number of such alternating green-red paths is at most $|M^*_i|$ (indeed, the path is completely determined by the first red edge it uses). 
    Therefore the number of edges of the reduced configuration which are exclusively blue is at least $C-|V(M^*_i)| - k |M^*_i| = C-(k+2)|M^*_i| \ge C - (k+2)km = (m+3)k \ge (|Z|+2)k+1$, where we used our choice of $C$. 
    Therefore Lemma~\ref{lem:embed} on input $Z$, $x=\phi(a_{i+1})$, $y=\phi(b_{i+1})$ and colour $c_{i+1}$ gives a red matching $M$ and we can set $M^*_{i+1}:=M^*_i \cup M$.
    (Observe that, by~\ref{eq:Q}, the neighbours of $x$ and $y$ in colour $c_{i+1}$ are outside $Z$.)
    It is easy to see that $M^*_{i+1}$ satisfies~\ref{eq:Q}. 
\endproof
\medskip

\subsection{Proof of \cref{thm:abs}} 
We can now show the main result of this section: Configurations which, for each colour, contain enough edges of that colour only, admit a red perfect matching.

\lateproof{Theorem~\ref{thm:abs}}
Let $C:=C(k)$ be a large enough constant for the arguments below to hold and let $m \in \{k,k+1\}$ be even. Let $\cL_1,\dots,\cL_s$ be a list of all possible $k$-configurations of order $m$, up to isomorphism. Note that $s=O_k(1)$. 
For every $j \in [s]$, by Lemma~\ref{lem:abs}, there exists a $k$-configuration $\cA_j$ such that $\cA_j$ is an absorber for $\cL_j$. We may assume that the $\cA_j$'s have pairwise disjoint vertex sets, and we let $\cA^*:=\cA_1\cup \dots \cup \cA_s$.
(Note $\cA^*$ is a $k$-configuration.)
Now, by Lemma~\ref{lem:find}, we can find a partial red matching $M^*_1$ in $\cM$ such that the reduced configuration $\cM_1$ contains $\cA^*$ as a subconfiguration, and we let $U$ be its vertex set. 
Note that, since $\cA^*$ is a configuration, $\cM_1$ has no edges between $U$ and $W:=V(\cM_1)\sm U$. In $\cM_1[W]$, we run the greedy algorithm to find a partial red matching $M^*_2$ that covers all but $m$ vertices. Let $\cL$ be the reduced configuration on these $m$ vertices. Then it must be isomorphic to $\cL_{j_0}$ for some $j_0\in [s]$. Now, for each $j\neq j_0$, there is a totally equalizing red matching for $\cA_j$, and there is a totally equalizing red matching for $\cL_{j_0}\cup \cA_{j_0}$ (and thus for $\cL \cup \cA_{j_0}$). Let $M^*_3$ be the union of all these red matchings. 
Hence, $M^*_1\cup M^*_2 \cup M^*_3$ is a totally equalizing matching of $\cM$, and we can easily choose the remaining red edges to find a red perfect matching for~$\cM$.
\endproof

\section{Counting -- Proof of Theorem~\ref{thm:count}}
\label{sec:counting}

\cref{thm:count} was stated for edge-disjoint perfect matchings, but we can show that the conclusion still holds if the matchings overlap only mildly, which is crucial for the deduction of \Cref{thm:random}.

\begin{theorem}\label{thm:count robust}
    Let $k\in \bN$ be fixed.
    Then the following holds for all $n$ sufficiently large. Let $V$ be a set of $n$ vertices and suppose $M_1,\dots,M_k$ are perfect matchings on~$V$ such that $\sum_{i,j:i\neq j} |M_i\cap M_j| \le n^{1/8}$.
    Then the number of perfect matchings $M^*$ such that $M^*\cup M_i$ forms a Hamilton cycle for each $i\in [k]$ is $\Theta_k\left(n^{-k/2}(n/e)^{n/2}\right)$.
\end{theorem}

\cref{thm:count} is then an immediate corollary of \Cref{thm:count robust}.

As already explained in \cref{subsec:basic_concepts}, given a configuration, we construct the red matching iteratively and, as long as at least $k+2$ vertices are uncovered, we can choose the next red edge greedily.
In order to count the number of red perfect matchings, roughly speaking, we want to estimate the number of choices the greedy algorithm has in each step. It is much more convenient to count ordered matchings in this way, which simply differs from the unordered case by a factor of $(n/2)!$.
First, we give a rough sketch of how the counting proof works.

Let $\cM$ be the given configuration and $M^*=(e_1^*,e_2^*,\dots,e_{n/2}^*)$ be an ordered red perfect matching.
For $0 \le t \le n/2-1$, let $\cM_t$ be the reduced configuration of the matching $\Set{e_1^*,\dots,e_{t}^*}$ with the convention that $\cM_0:=\cM$.
Moreover, let $n_t:=n-2t$ denote the order of the reduced configuration at time~$t$, and let $M_{t,i}$ denote the edges of $\cM_t$ of colour $i$.
The number of possibilities for $e_{t+1}^*$ is then the number of edges in the reduced configuration that are not occupied by any reduced matching, so it is exactly
\begin{equation*}
    \binom{n_t}{2}-|\cup_{i\in[k]} M_{t,i}|\, .
\end{equation*}
By assuming for now that the matchings stay (nearly) edge-disjoint, this amounts roughly to $$\binom{n_t}{2}-k\frac{n_t}{2}= \frac{n_t^2}{2} \cdot \left[ 1-\frac{k+1}{n_t}\right] \approx \frac{n_t^2}{2} \cdot \exp\left(-\frac{k+1}{n_t}\right) \,,$$ 
where we used that $1-x \approx \exp\left(-x\right)$.
By multiplying the choices we get
\begin{align}
\label{eq:counting}
\prod_{t=0}^{n/2-1} \left[ \frac{n_t^2}{2} \cdot \exp\left(-\frac{k+1}{n_t}\right) \right] =\left( \prod_{t=0}^{n/2-1} \frac{n_t^2}{2}\right) \cdot \exp\left(-(k+1)\sum_{t=0}^{n/2-1}\frac{1}{n_t}\right) \, .
\end{align}
We can easily see that
\begin{equation*}
    \prod_{t=0}^{n/2-1} \frac{n_t^2}{2} = \frac{(2^{n/2} \cdot (n/2)!)^2}{2^{n/2}} = 2^{n/2} \cdot (n/2)!^2 \, .
\end{equation*}
Moreover, $\sum_{t=0}^{n/2-1}1/n_t$ is the harmonic sum but only over the even integers, and thus it is approximately $\log (n/2)/2 \approx \log(n)/2$ (cf.~\cref{fact:series_1/n}) and hence the second term in~\eqref{eq:counting} is approximately $n^{-(k+1)/2}$.
Putting all together, the number of ordered red matchings is approximately $2^{n/2} \cdot (n/2)!^2 \cdot n^{-(k+1)/2}$.
After dividing this quantity by $(n/2)!$, we get that the number of unordered red matchings is approximately
\[
    2^{n/2} \cdot (n/2)! \cdot n^{-(k+1)/2} \approx n^{-k/2}(n/e)^{n/2}\, ,
\]
where we used Stirling's approximation $n! \approx n^{1/2} \cdot (n/e)^{n}$.

We will argue that this heuristic argument can be made rigorous. 
For that, there are two main obstacles arising from the assumptions we made in the simplified argument above.
Firstly, we cannot assume that the matchings stay edge-disjoint:
In fact, the number of choices at each step heavily depends on the already chosen red edges.
Secondly, the greedy algorithm does not necessarily produce a red perfect matching, as it is only guaranteed to run as long as at least $k+2$ vertices are left.
The details are different for the lower and upper bound.

For the lower bound, overlap actually helps us in some sense, since then the number of choices for the next red edge is larger and it will be enough to use the bound $|\cup_{i\in[k]} M_{t,i}| \le k \frac{n_t}{2}$. So the main concern is that towards the end of the process, we do not get stuck. That is, we want to stop the process at some time $m$, where $n/2-m=\Theta_k(1)$, and be left with a reduced configuration that has small enough overlap so that we can complete the red matching into a red perfect matching by an application of Theorem~\ref{thm:abs}. While this is not always true, we show that it is true in expectation.
Therefore, by Markov's inequality, a constant fraction of the outcomes can indeed be completed.

For the upper bound, we can ignore the problem that some of the red matchings of size $m$ might not be completable to a red perfect matching. 
However, we do have to control carefully the total contribution of the term $|\cup_{i\in[k]} M_{t,i}|$ over time.
Again, we do so for a uniformly random ordered red matching, where the overlap at time $t$ corresponds to the overlap of the reduced configuration of the first $t$ edges of the ordered red matching.
Intuitively, this is similar to running the greedy process at random, meaning that at each step we choose the next red edge uniformly at random among all the available ones.
However, they do not have the same distribution and we need the following statistical result.

\begin{theorem}\label{thm:stats}
    Let $k \in \mathbb{N}$ be fixed.
    Then the following holds for $n$ sufficiently large.
    Let $\cM$ be a $k$-configuration of order $n$ with matchings $M_1, \dots, M_k$ and suppose that $\sum_{i,j:i\neq j}|M_i \cap M_j|  \le n^{1/8}$. Let $M^*$ be a red perfect matching chosen uniformly at random. Then, for every edge $e\notin E(\cM)$, we have $$\prob{e \in M^*}=\Theta_k(1/n).$$
\end{theorem}

This is actually a corollary of Theorem~\ref{thm:count robust}. 
Indeed, from \Cref{thm:count robust} we know that the number of \emph{all} red perfect matchings of $\cM$ is $\Theta_k\left(n^{-k/2}(n/e)^{n/2}\right)$.
Then consider the reduced configuration $\cM'$ obtained by taking a fixed $e \not\in E(\cM)$ as red edge. 
By Item~\ref{prop_P_2} of \Cref{prop:counting_double_edges}, the total overlap increases by at most $2k^2$ and thus we can apply \Cref{thm:count robust} again and get that the number of red perfect matchings of $\cM'$ is $\Theta_k\left(n^{-k/2}(n/e)^{n/2-1}\right)$.
Therefore we can conclude that
$$ 
   \prob{e \in M^*}=\frac{\Theta_k\left(n^{-k/2}(n/e)^{n/2-1}\right)}{\Theta_k\left(n^{-k/2}(n/e)^{n/2}\right)} = \Theta_k(1/n)\, .
$$
In fact, we will not prove Theorem~\ref{thm:stats} directly since this seems roughly as hard as Theorem~\ref{thm:count robust}. Instead, we will prove the following variant for almost-perfect matchings, which suffices for our purposes.

\begin{lemma}\label{lem:stats_unordered}
Let $1/n \ll 1/D \ll 1/k$ with $k,D,n \in \mathbb{N}$. 
Let $\cM$ be a $k$-configuration of order $n$ and $M^*$ be a red matching of size $m:=n/2-D$ chosen uniformly at random. Then, for every edge $e\notin E(\cM)$, we have $$\prob{e\in M^*} = \frac{1 \pm o_D(1)}{n}\, .$$
\end{lemma}
The proof is postponed to \Cref{sec:switiching} and utilises the switching method.
We remark that for our application, it would be enough to have $\prob{e\in M^*} = \Theta_k(1/n)$, with the caveat of changing the constant $1/8$ in the bound of the statement of \Cref{thm:count robust} to an implicit one.
However, the proof of the weaker estimate via the same methods would not be simpler.

\subsection{General observations}
As already pointed out, in our calculation, we will encounter partial sums of the harmonic series and we will use the following simple fact. 
\begin{fact}
\label{fact:series_1/n}
    Let $n$ be a positive even integer and $s < n/2$ be an integer.
    Set $n_t:=n-2t$ for each $t \le s$.
    Then $\sum_{t=0}^{s-1} \frac{1}{n_t} = \frac{1}{2} \left(\log \frac{n}{n-2s}\pm O(1) \right)$.
\end{fact}

\proof
    Using the well-known fact that $\sum_{i=1}^n \frac{1}{i} = \log n \pm O(1)$, together with the observation that $n_t=2(n/2-t)$ and $n/2 -t \in \mathbb{N}$, we get that
    \begin{align*}
        \sum_{t=0}^{s-1} \frac{1}{n_t} &= \frac{1}{2}\left(\frac{1}{n/2}+\frac{1}{n/2-1}+\dots+\frac{1}{n/2-s+1}\right) \\
        &= \frac{1}{2} \left( \sum_{t=1}^{n/2} \frac{1}{t} -\sum_{t=1}^{n/2-s} \frac{1}{t} \right)
        = \frac{1}{2} \left(\log \frac{n/2}{n/2-s} \pm O(1)\right) = \frac{1}{2} \left(\log \frac{n}{n-2s} \pm O(1)\right)\, . 
    \end{align*}
\endproof

The next lemma provides the relevant estimates concerning the total overlap of a uniformly random ordered red matching. 
We recall that the total overlap of a configuration $\cM$ is denoted by $R(\cM)$ (cf.~\Cref{def:overlap}).

\begin{lemma}
\label{lem:sequence_has_nice_properties}
    Let $1/n \ll 1/D' \ll 1/D \ll 1/k$ with $k,D,n \in \mathbb{N}$, and fix a $k$-configuration $\cM$ on a vertex set of size $n$ with $R(\cM) \le n^{1/8}$.
    For an ordered red matching $M:=(e_1, \dots, e_{m})$ of size $m:=n/2-D$ and $t \in \{0,1,\dots,m\}$, set $n_t:=n-2t$ and denote by $R_t(M)$ the total overlap of the reduced configuration of the red matching $\{e_1,\dots,e_t\}$.
    Then, for a uniformly random ordered red matching $M$ of size $m$, we have $\mathrm{\mathbb{E}}\left[\sum_{t=0}^{m} \frac{R_t(M)}{n_t^2}\right] = O_{D'}(1)$, and moreover, for $m':=n/2-D'$, we have $\expn{R_{m'}(M)} < D/4$. 
\end{lemma}

In order to establish \Cref{lem:sequence_has_nice_properties}, we use the following result.

\begin{lemma}
\label{lem:sequence_R_t}
    Let $n \in \mathbb{N}$ and $m < n/2$ be an integer.
    Set $n_0:=n$ and, for each $t \in [m]$, set $n_t:=n-2t$.
    Let $R_0, R_1, \dots, R_{m}$ be random variables such that
    \begin{enumerate}[label=\rm{(S\arabic*)}]
        \item \label{prop_S_1} $\expn{R_{t}-R_{t-1}} \le k^2$ for each $t$; and
        \item \label{prop_S_2} there exist $0 < \alpha <2$ and $B>0$ such that $\expn{R_{t}} \le (1-\alpha/n_{t-1}) \cdot \expn{R_{t-1}}$ for each $t$ with $\expn{R_{t-1}} \ge B$.
    \end{enumerate}
    Then, for each $t \in [m]$, we have $\expn{R_t} = O_{\alpha}(1) \cdot \expn{R_0} \cdot \left(\frac{n_t}{n}\right)^{\alpha/2}+B+k^2$. Moreover $\expn{\sum_{t=0}^{m} \frac{R_t}{n_t^2}} = O_{\alpha}\left(\frac{\expn{R_0}}{n^{\alpha/2}}+B+k^2\right)$.
\end{lemma}

\proof   
    We start by the following claim.
    \begin{claim}
    \label{claim:exp_below_constant}
        If there exists $s \ge 0$ such that $\expn{R_{s+1}} \le B$, then $\expn{R_t} \le B+k^2$ for each $s+1 \le t \le m$.
    \end{claim}
    \claimproof
        This can be proved by induction on $t$. If $t=s+1$, the claim is true by assumption.
        Fix $t > s+1$.
        If $\expn{R_{t-1}} < B$, using~\ref{prop_S_1}, we have $\expn{R_t} \le \expn{R_{t-1}} + k^2 \le B+k^2$.
        If $\expn{R_{t-1}} \ge B$, then by~\ref{prop_S_2}, we have $\expn{R_t-R_{t-1}} \le 0$, and thus $\expn{R_t} \le \expn{R_{t-1}} \le B+k^2$.
    \endproof

    Let $s$ be the least $t$ with $0 \le t < m$ such that $\expn{R_{t+1}} \le B$ (if it exists, otherwise take $s=m$).
    Then, by Claim~\ref{claim:exp_below_constant}, $\expn{R_t} \le B+k^2$ for each $s+1 \le t \le m$.
    Therefore 
    \begin{equation}
    \label{eq:sum_R_t_1}
        \sum_{t=s+1}^{m} \frac{\expn{R_t}}{n_t^2} \le (B+k^2) \cdot \sum_{t=s+1}^{m} \frac{1}{n_t^2}= (B+k^2) \cdot O(1)\, ,
    \end{equation}
    where we used that $\sum_{i=1}^{\infty} i^{-2} < \infty$.
    
    On the other hand, by definition of $s$, for all $t \in [s]$, we have that $\expn{R_{t}} > B$ and thus
    \begin{align*}
        \expn{R_t} &\le \expn{R_0} \cdot \prod_{i=0}^{t-1} \left(1-\frac{\alpha}{n_i}\right) \le \expn{R_0} \cdot \exp\left(-\alpha \cdot \sum_{i=0}^{t-1} \frac{1}{n_i}\right) \\
        & \le
        \expn{R_0} \cdot \exp\left\{-\frac{\alpha}{2} \left(\log \frac{n}{n-2t} - O(1)\right)\right\} \\
        & =
         O_{\alpha}(1) \cdot \expn{R_0} \cdot \left(\frac{n-2t}{n}\right)^{\alpha/2}\, ,
    \end{align*}
    where the first inequality uses that from~\ref{prop_S_2} we have $\expn{R_{t}} \le \left(1-\frac{\alpha}{n_{t-1}}\right) \cdot \expn{R_{t-1}}$, the second inequality uses $1-x \le \exp(-x)$ and the third inequality uses \cref{fact:series_1/n}.
    The first part of the lemma follows from the calculation above together with Claim~\ref{claim:exp_below_constant}.
    
    Moreover, for all $t \in [s]$, we have
    \[
        \frac{\expn{R_t}}{n_t^2} = O_{\alpha}(1) \cdot \expn{R_0} \cdot \left(\frac{n-2t}{n}\right)^{\alpha/2} \cdot \frac{1}{(n-2t)^2} = O_{\alpha}(1) \cdot \frac{\expn{R_0}}{n^{\alpha/2}} \cdot \frac{1}{(n-2t)^{2-\alpha/2}}\, ,
    \]
    and thus 
    \begin{equation}
    \label{eq:sum_R_t_2}
        \sum_{t=0}^{s} \frac{\expn{R_t}}{n_t^2}= \frac{\expn{R_0}}{n^{\alpha/2}} \cdot O_{\alpha}(1)\, ,
    \end{equation}
    where we used that $\sum_{i=1}^{\infty} i^{-(2-\alpha/2)} < \infty$ as $2 - \alpha/2 >1$.

    The second part of the lemma then follows from~\eqref{eq:sum_R_t_1} and~\eqref{eq:sum_R_t_2}.
\endproof

We finish this subsection with the proof of \Cref{lem:sequence_has_nice_properties}, where we use that the total overlap of a uniformly random ordered red matching has a ``self-correcting'' behaviour in the following sense.
By \Cref{prop:counting_double_edges}, if the overlap is large, then there are many edges whose choice as the subsequent red edge would decrease the overlap.
Moreover, by \Cref{lem:stats_unordered}, every available edge is roughly equally likely to be contained in a uniformly random ordered red matching.
Therefore, the overlap decreases in expectation.

\lateproof{\Cref{lem:sequence_has_nice_properties}}
    Fix $t \in [m']$.
    By \cref{obs:random_stays_random}, conditioned on $e_1, \dots, e_{t-1}$, the ordered matching $(e_t,\dots,e_m)$ is distributed as a uniformly random ordered red matching of size $\tilde{m}:=m-(t-1)=n_{t-1}/2-D$ of the reduced configuration of $\{e_1,\dots, e_{t-1}\}$, which we denote by $\cM'$.
    The order of $\cM'$ is $n_{t-1} \ge n-2(m'-1) = 2D'+2$ and, since $1/D' \ll 1/D$, we can apply \cref{lem:stats_unordered} to $\cM'$.
    Therefore, for every edge $e \not\in E(\cM')$, we have 
    \begin{align*}
        \prob{e_t=e|e_1,\dots,e_{t-1}} &= \frac{\# \text{ ordered red matchings of $\cM'$ with first edge $e$}}{\# \text{ ordered red matchings of } \cM'} \\
        &=  \frac{(\tilde{m}-1)! \cdot \# \text{ unordered red matchings of $\cM'$ containing } e}{\tilde{m}! \cdot \# \text{ unordered red matchings of $\cM'$}} \\
        &= \frac{1}{\tilde{m}} \cdot \frac{1 \pm o_D(1)}{n_{t-1}}\\
        &=\frac{2 \pm o_D(1)}{n_{t-1}^2}\, ,
    \end{align*}
    where the third line uses \cref{lem:stats_unordered}.

    We now argue that the sequence of random variables $R_0(\cdot),\dots,R_{m'}(\cdot)$ satisfies the assumptions of \cref{lem:sequence_R_t} with $\alpha:=1/4$ and $B:=24k^2+2k^3$: Indeed, Property~\ref{prop_S_1} holds by Item~\ref{prop_P_2} of \cref{prop:counting_double_edges} (summing over all pairs of colours) and Property~\ref{prop_S_2} can be proved as follows.
    The calculation above shows that for every $t \in [m']$ and every available red edge $e$, 
    \begin{equation}
    \label{prop:Q_distribution_R_t}
        n^{-2}_{t-1} \le \prob{e_t=e|e_1,\dots,e_{t-1}} \le 3n_{t-1}^{-2}\, .
    \end{equation}
    For any distinct colours $i,j \in [k]$, set $R_{ij}$ (resp. $R_{ij}'$) to be the overlap in colours $i, j$ of the reduced configuration of the first $t-1$ (resp. $t$) red edges.
        Recall the definition of good and bad edges in \Cref{def:good_edges} and condition on $R_{t-1}$ (and thus on $R_{ij}$).
        From Property~\ref{prop_P_1} of \Cref{prop:counting_double_edges}, there are at least $(R_{ij}/2-k/2)n$ good edges and at most $n$ bad edges (for colours $i$ and $j$).
        Together with Property~\ref{prop_P_2} and~\eqref{prop:Q_distribution_R_t}, we have
        \begin{align*}
            \expn{R_{ij}'|R_{t-1}} & \le R_{ij} + 2 \cdot 3 n^{-2} \cdot (\# \text{ bad edges}) - 1 \cdot n^{-2} \cdot  (\# \text{ good edges}) \\
            & \le R_{ij} + \frac{6+k/2-R_{ij}/2}{n}\, .
        \end{align*}
        By summing over all pairs of colours, and using that $R_t=\sum_{1 \le i < j \le k} R_{ij}$ and $R_{t-1}=\sum_{1 \le i < j \le k} R_{ij}'$, we get
        \begin{align*}
            \expn{R_t|R_{t-1}} \le R_{t-1} + \frac{(6+k/2)k^2-R_{t-1}/2}{n}\, ,
        \end{align*}
        and taking the expectation (over $R_{t-1}$) on both sides gives
        \begin{align*}
            \expn{R_t} \le \expn{R_{t-1}} + \frac{(6+k/2)k^2-\expn{R_{t-1}}/2}{n} \le (1- \alpha/n) \cdot \expn{R_{t-1}}\, ,
        \end{align*}    
        where the last inequality follows from the choice of $\alpha$ and $B$ when assuming that $\expn{R_{t-1}}\ge B$.
    Then \cref{lem:sequence_R_t} implies that
    \[
    \expn{R_{m'}(M)} = O(1) \cdot R_0 \cdot \left(\frac{2D}{n}\right)^{1/8}+24k^2+2k^3+k^2 < D/4 \, ,
    \]
    and
    \[ 
        \mathrm{\mathbb{E}}\left[\sum_{t=0}^{m'} \frac{R_t(M)}{n_t^2}\right] = O\left(\frac{R_0}{n^{1/8}}+24k^2+2k^3+k^2\right) = O_{k}(1)\, , 
    \] 
    where we used that $R_0=R(\cM) \le n^{1/8}$ and $1/D \ll 1/k$.

    Since $R_t(M) \le \binom{k}{2} \frac{n_t}{2}$, we have $\mathrm{\mathbb{E}}\left[\frac{R_t(M)}{n_t^2}\right] \le k^2/4$ for any $t \in [m]$.
    Thus, by splitting the sum at $t=m'$, we can conclude that
    \[ 
        \mathrm{\mathbb{E}}\left[\sum_{t=0}^{m} \frac{R_t(M)}{n_t^2}\right] \le \mathrm{\mathbb{E}}\left[\sum_{t=0}^{m'} \frac{R_t(M)}{n_t^2}\right] + (m-m')k^2/4 = O_{D'}(1) \, . 
    \]
\endproof

\subsection{Lower bound}

\lateproof{the lower bound of \cref{thm:count robust}}
Choose new constants $D$ and $D'$ such that $$1/n \ll 1/D' \ll 1/D \ll 1/k\, ,$$
and set $C:=D/2$.
Let $\cM_0:=\cM$ be the given configuration, $m:=n/2-D$, $n_0:=n$ and, for $t \in [m]$, $n_t:=n-2t$.

Let $\cX$ be the set of ordered red matchings $(e_1^*,\dots,e_m^*)$ of size $m$.
We can construct the elements of $\cX$ greedily as follows:
Having chosen $e_1^*, \dots, e_{t}^*$, we choose $e_{t+1}^*$ to be an available red edge of the reduced configuration of $e_1^*, \dots, e_{t}^*$ (which has $n_t$ vertices).
Observe that different outcomes give different ordered matchings.
Moreover, the number of choices for $e_{t+1}^*$ is at least $\binom{n_t}{2}-k \frac{n_t}{2} = \frac{n_t^2}{2} \left(1-\frac{k+1}{n_t}\right) \ge \frac{n_t^2}{2} \cdot \exp\left(-\frac{k+1}{n_t}-2\frac{(k+1)^2}{n_t^2}\right)$, where we used that $1-x \ge \exp(-x-2x^2)$ for each $x \in [0,1/4]$ (which is satisfied in our case as $(k+1)/n_t \le (k+1)/(2D) \ll 1$ for $t \in [m]$).
Therefore we have 
\begin{align*}
    |\cX| &\ge \prod_{t=0}^{m-1} \left[\binom{n_t}{2}-k \frac{n_t}{2}\right] \\
    &= \Omega_D(1) \cdot \prod_{t=0}^{n/2-1} \frac{n_t^2}{2} \cdot \exp\left(-\frac{k+1}{n_t}-2\frac{(k+1)^2}{n_t^2}\right) \\
    &= \Omega_D(1) \cdot \left(\prod_{t=0}^{n/2-1} \frac{n_t^2}{2}\right) \cdot \exp\left(-(k+1) \sum_{t=0}^{n/2-1} \frac{1}{n_t}\right) \cdot \exp\left(-2(k+1)^2 \sum_{t=0}^{n/2-1} \frac{1}{n_t^2}\right) \, ,
\end{align*}
where the second line uses the inequality above and that if $m < t < n/2$, then $2 \le n_t \le 2D$, while the third line is just rearranging.
Using \cref{fact:series_1/n}, together with $\prod_{t=0}^{n/2-1} \frac{n_t^2}{2}=2^{n/2} \cdot (n/2)!^2$ and $\sum_{i=1}^{\infty} i^{-2} < \infty$, we have 
\begin{align*}
    |\cX| &= \Omega_D(1) \cdot 2^{n/2} \cdot (n/2)!^2 \cdot \exp\left(-(k+1) \frac{\log(n/2)-O(1)}{2}\right) \\
    & = \Omega_D(1) \cdot 2^{n/2} \cdot (n/2)!^2 \cdot n^{-(k+1)/2}\, .
\end{align*}

For a matching of $\cX$, we define its \defn{cut matching} to be the ordered red matching obtained by considering only the first $n/2-D'$ edges.
We now estimate the proportion of matchings in $\cX$ such that their cut matching can be completed to a full red matching.
Let $M_1^* \in \cX$ and $\tilde{\cM}:=\{\tilde{M_1},\dots,\tilde{M_k}\}$ be the reduced configuration of its cut matching $\tilde{M_1^*}$, say on vertex set~$\tilde V$ (which has size $2D'$).
Suppose that
\begin{equation}
\label{eq:matching_can_be_completed}
    \left|\tilde{M_i} \setminus \bigcup_{j \neq i} \tilde{M_j}\right| \ge C \text{ for each } i \in [k]\, .
\end{equation}
Then, by \cref{thm:abs}, $\tilde \cM$ has a red perfect matching $M_2^*$ (covering $\tilde V$), which we arbitrarily order.
By concatenating (in this order) $\tilde{M_1^*}$ and $M_2^*$, we obtain an ordered red perfect matching $M^*$ and note that different $\tilde{M_1^*}$ give different $M^*$.
Also observe that if there exists $i \in [k]$ such that $\left|\tilde{M_i} \setminus \bigcup_{j \neq i} \tilde{M_j}\right| < C$, then $R(\tilde{\cM}) \ge \sum_{j:j \neq i} |\tilde{M_j} \cap \tilde{M_i}| \ge D-C >0$.

By \Cref{lem:sequence_has_nice_properties}, for a uniformly random ordered red matching of size $m$, we have that the reduced configuration at time $n/2-D'$ (which we denoted by $\tilde \cM$) satisfies $\expn{R(\tilde \cM)} < D/4$.
By Markov's inequality, 
\[ 
    \prob{\tilde \cM \text{ does not satisfy }~\eqref{eq:matching_can_be_completed}} \le \prob{R(\tilde \cM) \ge D-C} \le \frac{\expn{R(\tilde \cM)}}{D-C} < \frac{1}{2}\, ,
\]
where we also used that $C=D/2$.

Therefore at least half of the partial ordered red matchings in $\cX$ have the property that their cut matching can be completed to a full ordered red perfect matching.
Moreover, given any ordered red matching of size $n/2-D'$, there are at most $O_{D'}(1)$ ordered red matchings in $\cX$ having it as their cut matching.  
Recalling that these matchings are ordered, we conclude that the number of unordered red perfect matchings is at least $\Omega_{D'}(1) \cdot \frac{|\cX|/2}{(n/2)!}= \Omega_{D'}(1) \cdot 2^{n/2-1} \cdot (n/2)! \cdot n^{-(k+1)/2}$.
This finishes the proof of the lower bound of \cref{thm:count robust} as $n!=\Theta(n^{1/2}\cdot (n/e)^{n})$ by Stirling's approximation.
\endproof

\subsection{Upper bound}

As already mentioned, a common technique to prove upper bounds on the number of certain combinatorial objects is the entropy method. Usually it is carried out from the outset for each individual problem. We make an attempt to formulate a general lemma that could potentially be applied to other problems as well. Although the proof steps are just the usual entropy arguments, we find the general statement quite appealing and, to the best of our knowledge, it has not been stated explicitly before. 

The setup is the following. 
Let $\Omega_1, \dots, \Omega_m$ be finite sets and $\cX \subseteq \Omega_1 \times \dots \times \Omega_m$. Our goal is to give an upper bound on the size of $\cX$. A natural approach is to count the number of possibilities for $x=(x_1,\dots,x_m)\in \cX$ by choosing $x_1, x_2, \dots, x_m$ sequentially one after the other and accounting for the number of choices available at each step. The issue is that having chosen $x_1, \dots, x_i$, the number of choices for $x_{i+1}$ can heavily depend on the previous choices.
A bound on $|\cX|$ can still be obtained by considering the worst case, but this provides bad estimates. Our goal is to consider instead the average case.

One issue with formalizing a general statement is that the definition of a suitable choice in a specific step depends heavily on the concrete problem. 
To overcome this, we make the following abstract definition, which defines a suitable choice simply as a choice for which there is at least one completion to a full element of~$\cX$.

For $x:=(x_1,\dots,x_m)\in \cX$ as above and $t \in [m]$, define 
\[
    \cS_t(x):=\left\{z \in \Omega_{t}:     
    \begin{array}{c}
    \text{there exists } x' \in \cX \text{ with } x'_{t}=z \text{ and}\\
    x'_{i}=x_{i} \text{ for each } i \in [t-1]
    \end{array}
    \right\}\, ,
\] 
and observe that $\cS_t(x) \neq \emptyset$ (since $x_{t} \in \cS_t(x)$).
Moreover $\cS_t(x)$ actually depends only on $x_{1}, \dots, x_{t-1}$ and thus we will denote it by $\cS(x_{1}, \dots, x_{t-1})$ as well (and by $\cS_1$ when $t=1$).
Finally, define $s_t(x):=|\cS_t(x)|$. Trivially, we have $|\cX| \le \prod_{t \in [m]} \max_{x \in \cX} \left[  s_t(x) \right]$.
The following lemma now says that instead of considering the worst case in each step, we can take the average over a random element from~$\cX$. We remark that, although the definition of $s_t(x)$ is quite abstract, in practice one can easily obtain upper bounds on $s_t(x)$ by discarding all choices that are obviously not completable due to some condition that is specific to the problem at hand.

\begin{lemma}
\label{lem:entropy_counting}
    Let $\Omega_1, \dots, \Omega_m$ be finite sets and $\cX \subseteq \Omega_1 \times \dots \times \Omega_m$ with $\cX \neq \emptyset$.
    Then
    \[
        |\cX| \le \prod_{t \in [m]} \mathrm{\mathbb{E}}_{x \in \cX}\left[  s_t(x) \right]\, .
    \]
\end{lemma}

Note that in the proof we do not explicitly use the notion of entropy. We only rely on the AM-GM inequality and basic laws of probability, but the proof steps naturally correspond to the usual entropy proofs.

\lateproof{Lemma~\ref{lem:entropy_counting}}
Define a random variable $Z = (z_1,\dots,z_m) \in \cX$ as follows.
We let $z_{1}$ be a uniformly random element from~$\cS_1$. 
(That is, $z_{1}$ is an element of $\Omega_{1}$ chosen uniformly at random among those which can be completed to an element of $\cX$.)
Then, assuming we have chosen $z_{1}, \dots, z_{t-1}$, we let $z_{t}$ be a uniformly random element from $\cS(z_{1}, \dots, z_{t-1})$.
Observe that $Z \in \cX$ by construction.

    For $x \in \cX$, we have
    \begin{align*}
        \prob{Z=x} = \prob{\bigcap_{t \in [m]} \{z_{t}=x_{t}\}}
        & = \prod_{t \in [m]} \prob{z_{t}=x_{t} \Big| 
        \begin{array}{c}
            z_i=x_{i} \\
            \text{ for each } i \in [t-1]
        \end{array}}\\
        &= \prod_{t \in [m]} \frac{1}{s_t(x)}\, , 
    \end{align*}
    where we recall that since $x \in \cX$, we have $x_{t} \in \cS_t(x)$ and thus $s_t(x) \neq 0$ for all $t$.
    Therefore,
    \begin{align*}
        \log|\cX| &= \log \frac{|\cX|}{\sum_{x \in \cX} \prob{Z=x}} \le \frac{1}{|\cX|} \cdot \log \frac{1}{\prod_{x \in \cX} \prob{Z=x}} \\
        &= \frac{1}{|\cX|} \cdot \log\left( \prod_{x \in \cX} \prod_{t \in [m]} s_t(x)\right) = \sum_{t \in [m]} \frac{1}{|\cX|} \sum_{x \in \cX} \log\left( s_t(x)\right) \le \sum_{t \in [m]} \log\left( \frac{1}{|\cX|} \sum_{x \in \cX} s_t(x)\right) \, ,
    \end{align*}
    where the first equality uses that, since $Z \in \cX$, we have $\sum_{x \in \cX} \prob{Z=x} = 1$, the first inequality uses the AM-GM inequality, and the second inequality follows from Jensen's inequality and the concavity of the $\log$-function.
    By taking the exponential on both sides, we get
    \[
        |\cX| \le \prod_{t \in [m]} \frac{1}{|\cX|} \sum_{x \in \cX} s_t(x) = \prod_{t \in [m]} \mathrm{\mathbb{E}}_{x \in \cX}\left[  s_t(x) \right]\, ,
    \]
    as desired.
\endproof

\lateproof{the upper bound of \cref{thm:count robust}}
    Choose new constants $D$ and $D'$ such that $$1/n \ll 1/D' \ll 1/D \ll 1/k\, .$$
    Set $m:=n/2-D$, $n_0:=n$ and, for $t \in [m]$, $n_t:=n-2t$.
    
    Let $\cX$ be the set of ordered red matchings of size $m$.
    For an ordered matching $M:=(e_1, \dots, e_{m}) \in \cX$ and $t \in [m]$, define $s_t(M)$ to be the number of edges $e$ of $K_n$ such that there exists an ordered matching $M':=(e_1', \dots, e_{m}') \in \cX$ with $e_i'=e_i$ for each $i \in [t-1]$ and $e_{t}'=e$. 
    Moreover, denote by $\cM_{t-1}(M)$ the reduced configuration of the red matching $\{e_1, \dots, e_{t-1}\}$.
    Then, we have
    \[
        |\cX| \le \prod_{t \in [m]} \mathrm{\mathbb{E}}_{M \in \cX}\left[ s_t(M) \right] \le \prod_{t=0}^{m-1} \mathrm{\mathbb{E}}_{M \in \cX}\left[ \binom{n_t}{2}- e(\cM_t(M)) \right]\, ,
    \]
    where the first inequality uses \cref{lem:entropy_counting} and the second one follows from the fact that if an edge is counted by $s_{t}(M)$, then it must be an available edge of $\cM_{t-1}(M)$.
    For a more pleasant notation, we drop the subscript in $\mathrm{\mathbb{E}}_{M \in \cX}[\cdot]$ from now on, since we will always take expectations with respect to a uniformly random matching $M$ from~$\cX$.
    
    Letting $\left(\cM_t(M)\right)_i$ be the set of edges of colour $i$ of $\cM_t(M)$ and setting $R_t(M):=R(\cM_t(M))$, note that 
    \begin{align*}
        e(\cM_t(M)) &=
        |\cup_{i \in [k]} (\cM_t(M))_i| \\
        &\ge \sum_{i \in [k]} |(\cM_t(M))_i| - \sum_{1 \le i <j\le k} |(\cM_t(M))_i \cap (\cM_t(M))_j|=k \frac{n_t}{2}- R_t(M)\, ,
    \end{align*}
    and thus $\binom{n_t}{2}- e(\cM_t(M)) \le \binom{n_t}{2}-k \frac{n_t}{2}+ R_t(M) = \frac{n_t^2}{2} \cdot \left[1-\frac{k+1}{n_t}+\frac{R_t(M)}{n_t^2}\right]$.   
    Then
    \begin{align*}
     |\cX| & \le \prod_{t=0}^{m-1} \frac{n_t^2}{2} \cdot \left[1-\frac{k+1}{n_t}+\frac{\expn{R_t(M)}}{n_t^2}\right] \\
     & \le \prod_{t=0}^{m-1} \frac{n_t^2}{2} \cdot \exp\left(-\frac{k+1}{n_t}+\frac{\expn{R_t(M)}}{n_t^2}\right) \\
     &= \left(\prod_{t=0}^{m-1} \frac{n_t^2}{2}\right) \cdot \exp\left(-(k+1) \sum_{t=0}^{m-1} \frac{1}{n_t}\right) \cdot \exp\left(\expn{\sum_{t=0}^{m-1} \frac{R_t(M)}{n_t^2}}\right) \, ,
    \end{align*}
    where we used that $1+x \le e^x$ for each $x \in \mathbb{R}$ in the second line.
    Observe that $\prod_{t=0}^{m-1} \frac{n_t^2}{2} \le \prod_{t=0}^{n/2-1} \frac{n_t^2}{2} = 2^{n/2} \cdot (n/2)!^2$ and, by \cref{fact:series_1/n}, it holds that $\exp\left(-(k+1) \sum_{t=0}^{m-1} \frac{1}{n_t}\right) =O_D(1) \cdot n^{-(k+1)/2}$.
    Moreover, by \cref{lem:sequence_has_nice_properties}, $\mathrm{\mathbb{E}}\left[\sum_{t=0}^{m} \frac{R_t(M)}{n_t^2}\right] = O_{D'}(1)$.
    
    We can then conclude that $|\cX| = O_{D'}(1) \cdot 2^{n/2} \cdot (n/2)!^2 \cdot n^{-(k+1)/2}$.  
    Finally observe that if $M:=(e_1,\dots, e_{n/2})$ is an ordered red perfect matching, then $(e_1,\dots,e_{m})$ is an ordered red matching of size $m$. Moreover, every ordered matching of size $m$ can be completed to an ordered red perfect matching in $O_D(1)$ ways.
    Therefore the number of ordered red perfect matchings is $O_D(1) \cdot |\cX|$.
    After dividing this bound by $(n/2)!$ to account for the ordering, we conclude that the number of unordered red perfect matchings is $O_{D'}(1) \cdot 2^{n/2} \cdot (n/2)! \cdot n^{-(k+1)/2}$.
    This finishes the proof of the upper bound of \cref{thm:count robust} as $n!=\Theta(n^{1/2}\cdot (n/e)^{n})$ by Stirling's approximation.
\endproof

\section{Switching lemma}
\label{sec:switiching}
In this section, we prove \cref{lem:stats_unordered} by using the \defn{switching method}, which relies on the following strategy.
Suppose we have two disjoint finite sets $S, T$, and an operation, usually called switch, that takes an object in $S$ and modifies it to make an object in $T$.
If $d_S$ is the average number of switchings that can be applied to an object in $S$ and $d_T$ is the average number of switchings that can be applied to an object in $T$, then $d_S|S| = d_T |T|$. 
An estimate of the relative values of $d_S$ and
$d_T$ provides an estimate of the relative sizes of $S$ and~$T$. This is very useful in particular when it is hard to determine the total sizes of $S$ and $T$ directly.

We now move to our setting.
For the rest of the section, we fix a $k$-configuration $\cM$ of order $n$ and an edge $e:=ab \not\in E(\cM)$.
We assume that
\[
    1/n \ll 1/D \ll 1/k\, ,
\]
denote by $\cR$ the set of red matchings of $\cM$ of size $m:=n/2-D$, and define the following partition of $\cR$:
\begin{equation}
\label{eq:partition_classes}
\begin{gathered}
    A:=\{M \in \cR: e \in M\}\, ,  \\
    X:=\{M \in \cR: e \not\in M \text{ and } a, b \in V(M)\}\, , \\        
    Y:=\{M \in \cR: e \not\in M, \text{ and either } a \text{ or } b \in V(M)\}\, ,\\
    Z:=\{M \in \cR: a, b \not\in V(M)\}\, .
\end{gathered}
\end{equation}

We first prove the upper bound of \Cref{lem:stats_unordered}, which is relatively straightforward and serves to illustrate the method.

\lateproof{the upper bound of \Cref{lem:stats_unordered}}
It is enough to compare $|A|$ and $|X|$.
Define the bipartite graph with parts $A$ and $X$ as follows:
For $M \in A$ and $M' \in X$, we add the edge $MM'$ if $M'$ can be obtained from $M$ by deleting $e$ and another red edge, and adding two new red edges incident to a vertex of $V \setminus V(M)$ and to $a$ and $b$, respectively. 
Then, for every $M \in A$, we have $\deg_H(M) \ge (n/2-D-1)(2D-k)(2D-k-1) \ge (1-o_D(1))2D^2n$, as we must remove $e$ (this is a fixed choice) and another red edge (for which we have $n/2-D-1$ choices), then add an
available red edge from $a$ to a vertex in $V \setminus V(M)$ (for which we have at least $2D-k$ choices) and, after that, an available red edge from $b$ to a vertex in $V \setminus V(M)$ not yet used for $a$ (for which we have at least $2D-k-1$ choices).
Similarly, let $M' \in X$ and $a', b' \in V(M')$ such that $aa',bb' \in M'$. We have $\deg_H(M') \le \binom{2D+2}{2}$, as the edges to be removed are fixed and we must add $e$ (if possible) and a red edge in $(V \setminus V(M')) \cup \{a',b'\}$ (for which we have at most $\binom{2D+2}{2}$ choices).
Then $|A| \cdot \min_{M \in A} \deg_H(M) \le e(H) \le |X| \cdot \max_{M' \in X} \deg_H(M')$ and thus 
\begin{equation}
\label{eq:upper_bound_A}
    \frac{|A|}{|X|} \le \frac{\binom{2D+2}{2}}{(1-o_D(1))2D^2n} = \frac{1 + o_D(1)}{n}\, .
\end{equation}

The result follows as $\prob{e \in M^*} = \frac{|A|}{|A|+|X|+|Y|+|Z|} \le \frac{|A|}{|X|} \le \frac{1 + o_D(1)}{n}$.
\endproof

Unfortunately, such an easy switching argument does not go through for the lower bound. Indeed, a lower bound on $|A|/|X|$ would require a lower bound on $\deg_H(M')$ (where we are using the same notation as in the proof above). 
The issue is that some $M' \in X$ may be isolated in the graph $H$ as, after removing the two red edges incident to $a$ and $b$, it is not guaranteed that we can add $e$ as a red edge.
The main work of this section is to analyse the proportion of matchings for which $e$ is ``blocked''. Since we anyways have to delete the red edges at $a$ and $b$ before we can possibly add $e$, we will only consider matchings in~$Z$ for this.

\begin{defin}[Blocked edge and blocking path]
    Given a colour $i$ and a red matching $M$ with $a,b \not\in V(M)$, we say that $e$ is \defn{$i$-blocked} (in $M$) if $e$ appears in colour $i$ in the reduced configuration of $M$, that is, if the maximal alternating $i$-red path starting at $a$ ends in $b$. 
    We call such path the \defn{$i$-blocking path}.
    We say that $e$ is \defn{blocked} if there exists a colour $i$ such that $e$ is $i$-blocked.
\end{defin}

\begin{figure}[htp]
    \centering
    \includegraphics[scale=0.8]{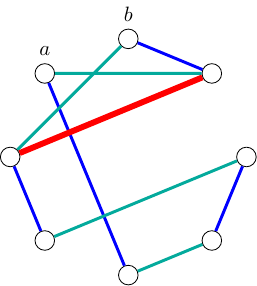}
    \caption{Let $e=ab$. Then $e$ is blocked in colour green, but not blocked in colour blue.}
\end{figure}

Note that for a matching $M\in Z$ and a colour $i$, the graph induced by $M$ and the matching of~$\cM$ of colour $i$ decomposes into a collection of vertex-disjoint alternating $i$-red paths, with each path starting and ending with an edge of colour~$i$.
(Observe that we may have paths consisting of a single edge of colour $i$.)
Moreover the endpoints of these paths are exactly the $2D$ vertices not covered by~$M$, and the path containing $a$ contains $b$ as well if and only if $e$ is $i$-blocked. Hence, heuristically, we would expect that for a random matching $M$ from $Z$, the probability that $e$ is $i$-blocked is $O(1/D)$.
We formalise this argument after introducing the relevant definitions.

For each colour $i \in [k]$, we define
\begin{equation*} 
\begin{gathered}
B^i:=\{M \in \cR: a, b \not\in V(M) \text{ and $e$ is $i$-blocked}\}\, ,\\
N^i:=\{M \in \cR: a, b \not\in V(M) \text{ and $e$ is not $i$-blocked}\}\, .
\end{gathered}
\end{equation*}
Then we let $B:=\cup_{i \in [k]} B^i$ and $N:= Z \setminus B$.
Note that a red matching $M \in \cR$ belongs to $N$ if and only if $a,b \not \in V(M)$ and $e$ is not blocked in any colour.
\begin{lemma}
\label{lem:B_and_N}
    For each colour $i \in [k]$, we have $|B^i| = O(1/D) \cdot |Z|$.
\end{lemma}

Fix a colour $i \in [k]$.
For an integer $\ell \ge 1$, we set $B_{\ell}^i$ to be the subset of $B^i$ consisting of the matchings for which the $i$-blocking path has exactly $\ell$ red edges.
Similarly, we let $N_{\ell}^i$ to be the subset of $N^i$ consisting of the matchings for which the maximal alternating $i$-red path starting at $a$ and the one starting at $b$ have in total exactly $\ell$ red edges.
(Recall these two paths are distinct by definition of $N^i$.)
The first step is to compare $B^i_\ell$ with $N^i_{\ell-1}$.

\begin{lemma}
\label{lem:B_ell_N_ell-1}
    For each colour $i \in [k]$ and integer $\ell \ge 1$, we have $|B_{\ell}^i| \le \frac{n}{\ell \cdot \Omega(D^2)} \cdot |N_{\ell-1}^i|$.
\end{lemma}

\begin{figure}[htp]
    \centering
    \includegraphics[scale=0.7]{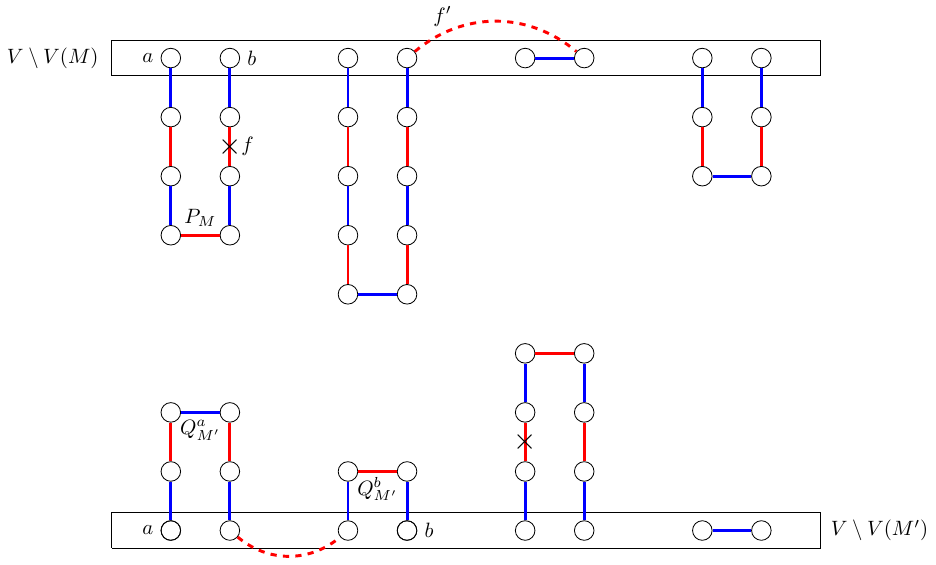}
    \caption{Switching between $B^i_\ell$ and $N^i_{\ell-1}$ in \cref{lem:B_ell_N_ell-1}.}
    \label{fig:first_lemma}
\end{figure}

\proof
    Suppose $i$ is the colour blue.
    For $M \in B_{\ell}^i$, let $P_M$ be the blue-blocking path.
    For $M' \in N_{\ell-1}^i$, let $Q_{M'}^a$ and $Q_{M'}^b$ be the maximal alternating blue-red paths starting at $a$ and $b$, respectively.
    Define $H$ to be the bipartite graph with parts $B_{\ell}^i$ and $N_{\ell-1}^i$ where $MM' \in E(H)$ if and only if $M'$ can be obtained from $M$ by removing a red edge $f$ of $P_M$ and adding a new red edge $f'$ with endpoints in $V \setminus (V(M) \cup \{a,b\})$ (cf.~\cref{fig:first_lemma}).
        
    Then, for every $M \in B^i_\ell$, we have $\deg_H(M)= \ell \cdot \Omega(D^2)$, as $f$ can be chosen in $\ell$ ways and $f'$ in at least $\binom{2D-2}{2}-k \cdot (D-1)$ ways. 
    Moreover, for every $M' \in N_{\ell-1}^i$, we have $\deg_H(M') \le n$ as, if possible, we must add the red edge joining the endpoints of $Q_{M'}^a$ and $Q_{M'}^b$ (different from $a$ and $b$), and remove a red edge of $M'$ not in the path $Q_{M'}^a$ or $Q_{M'}^b$ (for which we have at most $n$ choices). 
    
    The lemma then follows from $\ell \cdot \Omega(D^2) \cdot |B_{\ell}^i|=e(H) \le n \cdot |N_{\ell-1}^i|$.   
\endproof

When $\ell$ is linear (in $n$), \cref{lem:B_ell_N_ell-1} provides a good bound, but this is not the case anymore if $\ell$ is small.
However we are able to show that there are only few matchings in $B^i$ whose $i$-blocking path is short, and we do so by using the following lemma.

\begin{lemma}
\label{lem:prob_less_1/n}
    Let $\ell \le n/10$ and $T:=\{e_1, e_2, \dots, e_{\ell}\}$ be a red matching such that none of the edges $e_i$ is incident to $a$ or $b$.
    Then for a uniformly random red matching $M^*$ of size $n/2-D$ and any $x,y \not \in \cup_{i \in [\ell]} e_i$, it holds that
    \[
        \prob{xy \in M^* \Big| 
        \begin{array}{c} 
            a,b \not \in V(M^*) \text{ and }\\
            T \subseteq M^*
        \end{array}
        } = O\left( \frac{1}{n}\right)\, .
    \]
\end{lemma}

\proof
    The proof is identical to the upper bound of \Cref{lem:stats_unordered}, which we already proved at the beginning of the section.
    Let $\tilde{\cR}$ be the set of red matchings $M$ of size $n/2-D$ such that $a, b \not \in V(M)$ and $T \subseteq M$.
    Define $\tilde{A}:=\{M \in \tilde{\cR}: xy \in M\}$ and $\tilde{X}:=\{M \in \tilde{\cR}: xy \not\in M \text{ and } x, y \in V(M)\}$. 
    Then the required probability is at most $|\tilde{A}|/|\tilde{X}|$.
    
    Define the bipartite graph with parts $\tilde{A}$ and $\tilde{X}$ as follows: For $M \in \tilde{A}$ and $M' \in \tilde{X}$, we add the edge $MM'$ if $M'$ can be obtained from $M$ by deleting $xy$ and another red edge not in $T$, and adding two new red edges incident to a vertex of $V \setminus (V(M) \cup \{a,b\})$ and to $x$ and $y$, respectively. 
    Then, for every $M \in \tilde{A}$, we have $\deg_H(M) \ge (n/2-D-1-|T|)(2D-2-k)(2D-3-k)$. 
    Similarly, for every $M' \in \tilde{X}$, we have $\deg_H(M') \le \binom{2D}{2}$, as the edges to be removed are fixed and we must add $xy$ (if possible) and a red edges in $V \setminus (V(M') \cup \{a,b\})$. The result follows from $|\tilde{A}| \cdot \min_{M \in \tilde{A}} \deg_H(M) \le e(H) \le |\tilde{X}| \cdot \max_{M' \in \tilde{X}} \deg_H(M')$.    
\endproof

Finally we are able to prove \cref{lem:B_and_N}.

\lateproof{\cref{lem:B_and_N}}
    Let $M^*$ be a uniformly random red matching of size $n/2-D$ and set $\ell^*:=n/(10D)$.
    Partition $B^i=B^i_{< \ell^*} \cup B^i_{\ge \ell^*}$, where $B^i_{< \ell^*}:=\bigcup_{\ell < \ell^*} B^i_{\ell}$ and $B^i_{\ge \ell^*}:=\bigcup_{\ell \ge \ell^*} B^i_{\ell}$.

    Fix $\ell < \ell^*$ and a colour $i$.
    Define $\cE^i_\ell$ to be the event that there is an $i$-red alternating path starting with an edge of colour $i$ at $a$ and ending with an edge of colour $i$ at $b$ containing exactly $\ell$ red edges.
    We bound the probability of $\cE^i_\ell$ by conditioning on the first $\ell-1$ red edges on such a path.
    So let $\cT$ be the set of partial red matchings $T$ of size $\ell-1$ such that the maximal alternating $i$-red path starting with an edge of colour $i$ at $a$ does not contain $b$ and has exactly $\ell-1$ red edges (so these are precisely the edges of $T$).
    Let $a'$ be the endpoint of this path different from $a$, and $b'$ such that $bb'$ is an edge of colour $i$.
    Then, for each $T \in \cT$, we have 
    \[
        \prob{\cE^i_\ell \Big| 
        \begin{array}{c} 
            a,b \not \in V(M^*) \text{ and } \\
            T \subseteq M^*
        \end{array}
        } = 
        \prob{a' b' \in M^* \Big| 
        \begin{array}{c} 
            a,b \not \in V(M^*) \text{ and } \\
            T \subseteq M^*
        \end{array}} = O(1/n)\, ,
    \]
    where the last equality uses \cref{lem:prob_less_1/n}.
    Therefore, since if $\cE^i_\ell$ holds then $T \subseteq M^*$ for some $T \in \cT$ and since the events $\{T \subseteq M^*\}_{T \in \cT}$ are pairwise disjoint, by the law of total probability it holds that
    \[
        \prob{\cE^i_\ell | a,b \not \in V(M^*)} = \sum_{T \in \cT} \prob{\cE^i_\ell \Big| 
        \begin{array}{c} 
            a,b \not \in V(M^*) \text{ and } \\
            T \subseteq M^*
        \end{array}} \cdot \prob{T \subseteq M^*}
        = O(1/n)\, . 
    \]
    The last line is equivalent to $|B^i_{\ell}| = O(1/n) \cdot |Z|$, and this holds for each $\ell < \ell^*$.
    By a union bound over all $\ell < \ell^*$, we have $|B^i_{<\ell^*}| \le O(1/D) \cdot |Z|$.

    For $\ell \ge \ell^*$, we get from \cref{lem:B_ell_N_ell-1} that
    \[
        |B_{\ell}^i| \le \frac{n}{\ell \cdot \Omega(D^2)} \cdot |N_{\ell-1}^i| \le \frac{n}{\ell^* \cdot \Omega(D^2)} \cdot |N_{\ell-1}^i| = O(1/D) \cdot |N_{\ell-1}^i|\, .
    \]
    Therefore
    \[
        |B^i_{\ge \ell^*}| = \sum_{\ell \ge \ell^*} |B^i_\ell| = O(1/D) \cdot \sum_{\ell \ge \ell^*} |N^i_{\ell-1}| \le O(1/D) \cdot |Z|\, .
    \]  
    Overall we get $|B^i|=|B^i_{< \ell^*}| + |B^i_{\ge \ell^*}| = O(1/D) \cdot |Z|$.
\endproof

We are now ready to complete the proof of the switching lemma.
\lateproof{the lower bound of \cref{lem:stats_unordered}}
    We compare the sizes of the sets $A,X,Y,Z$, which partition $\cR$ and have been defined in~\eqref{eq:partition_classes}.
   
    We start with $A$ and $Z$, arguing via $N$.
    Define the bipartite graph $H$ with parts $A$ and $N$ as follows:
    For $M \in A$ and $M' \in N$, we add the edge $MM'$ if $M'$ can be obtained from $M$ by deleting $e$ and adding a red edge $f$ not incident to $a,b$.
    Then, for every $M \in A$, we have $\deg_H(M) \le \binom{2D}{2}$ as $f$ needs to be chosen in $V \setminus V(M)$. 
    Similarly, for every $M' \in N$, we have $\deg_H(M') \ge n/2-D$, as we must remove an edge of $M'$ and add $e$, which is always possible as $e$ is not blocked by definition of~$N$. 
    Then $|A| \cdot \binom{2D}{2} \ge e(H) \ge |N| \cdot (n/2-D)$ and thus $|A| \ge \frac{|N|}{4D^2} \cdot \left(1 -o_D(1)\right)\cdot n$.    
    By \Cref{lem:B_and_N} and a union bound over the $k$ colours, we have $|B|=O(k/D) \cdot |Z|$.
    Using $|Z|=|B|+|N|$, we conclude $|N|=(1-O(k/D))\cdot|Z|$.
    In particular, $|A| \ge \frac{|Z|}{4D^2} \cdot \left(1-o_D(1)\right) \cdot n$.

    Next, we compare the sizes of $X$ and $Z$.
    Define the bipartite graph $H$ with parts $X$ and $Z$ as follows: For $M \in X$ and $M' \in Z$, we add the edge $MM'$ if $M'$ can be obtained from $M$ by removing the two red edges incident to $a$ and $b$ and adding two red edges not incident to any of $a,b$. 
    Then, for every $M \in X$, we have $\deg_H(M) \ge \frac{1}{2} \cdot \left(\binom{2D+2}{2}-k(D+1)\right) \cdot \left(\binom{2D}{2}-kD\right)$
    and, for every $M' \in Z$, we have $\deg_H(M') \le \binom{n/2-D}{2} (2D+3)(2D+1)$.
    Arguing as above, we get $|X| \le \frac{|Z|}{4D^2} \cdot \left(1+o_D(1)\right) \cdot n^2$.      
    Together with the inequality~\eqref{eq:upper_bound_A}, this also implies that $|A| \le \frac{|Z|}{4D^2}\cdot (1+o_D(1)) \cdot n$.

    Finally, we compare the sizes of $Y$ and $Z$.
    Let $M \in Y$ and suppose without loss of generality that $a \in V(M)$ (and thus $b \not\in V(M)$).        
    Let $a' \in V(M)$ be such that $aa' \in M$.
    Define the bipartite graph $H$ with parts $Y$ and $Z$ as follows: 
    For $M \in Y$ and $M' \in Z$, we add the edge $MM'$ if $M'$ can be obtained from $M$ by removing $aa'$ and adding an edge with endpoints in $V \setminus (V(M) \cup \{b\})$.
    Then, for every $M \in Y$, we have $\deg_H(M) = \Omega_k(D^2)$. 
    Moreover, for every $M' \in Z$, we have $\deg_H(M') =O_D(n)$ as we must remove an edge of $M'$ (and we have $n/2-D$ choices), choose $a$ or $b$ (2 choices) and, assuming we chose $a$, add an available edge from $a$ to a vertex in $V \setminus (V(M') \cup \{b\})$ (for which the number of choices is at most $2D-2$). 
    Arguing as above, we get $|Y| \le |Z| \cdot O_k(n)$.

    Combining everything together, we have $|A|+|X|+|Y|+|Z| \le \frac{|Z|}{4D^2} \cdot (1+o_D(1)) \cdot n^2$, and thus
    \begin{align*}
        \prob{e \in M^*} = \frac{|A|}{|A|+|X|+|Y|+|Z|} \ge \frac{\frac{|Z|}{4D^2} \cdot (1-o_D(1)) \cdot n}{\frac{|Z|}{4D^2} \cdot (1+o_D(1)) \cdot n^2} \ge \frac{1-o_D(1)}{n}\, .
    \end{align*}
\endproof

\section{Spreadness and proof of Theorem~\ref{thm:random}}
\label{sec:proof_thm_random}
The proof of \cref{thm:random} uses the fractional version of the Kahn–Kalai conjecture~\cite{KK:07}.
This was proposed by Talagrand~\cite{talagrand:10} and recently resolved by Frankston, Kahn, Narayanan and Park~\cite{FKNP:21}. The Kahn–Kalai conjecture was then proved in full by Park and Pham~\cite{PP:23}, but the fractional version is sufficient for us.

Suppose one wants to show that w.h.p.~$G(n,p)$ contains a certain substructure.
Define an auxiliary hypergraph $\cH$, where the vertices of $\cH$ are the edges of $K_n$ and the hyperedges correspond to the desired substructures.
So the question is now whether a random subset of the vertices of $\cH$ contains some hyperedge.
The fractional Kahn-Kalai conjecture reduces this probabilistic problem to the `extremal' problem of showing that $\cH$ is spread, where spreadness is defined as follows. 
We use $|\cH|$ to denote the number of hyperedges of $\cH$.
\begin{defin}
Let $\cH$ be a hypergraph on $V$.
We say that $\cH$ is \defn{$q$-spread} if for every $S\In V$ we have that $\Big|\{H \in \cH: S \subseteq H\}\Big| \le q^{-|S|} \cdot |\cH|$, i.e.~if the number of hyperedges containing $S$ is at most a $q^{-|S|}$-fraction of all edges.
\end{defin}

The statement we use here follows from~\cite[Theorem 1.6]{FKNP:21} and is taken from~\cite[Theorem 1.2]{PSSS:22+}.

\begin{theorem}
\label{thm:spread}
There exists an absolute constant $C$ such that the following holds.
Let $\cH$ be an $N$-vertex $q$-spread hypergraph with $|\cH|>0$. Suppose that $W$ is a random subset of $V(\cH)$ where every vertex is included independently with probability $p\ge Cq^{-1} \log N$. Then w.h.p.~(as $N \to \infty$) there exists an edge of $\cH$ such that all its vertices have been selected.
\end{theorem}

We show how \cref{thm:spread} and \cref{thm:count robust} can be combined to prove \cref{thm:random}.

\lateproof{Theorem~\ref{thm:random}}
Let $C >0$ be large enough and $p \ge C (\log n)/n$.
We reveal $G(n,p)$ in $k$ rounds $G_1,\dots,G_k$, meaning that each $G_i$ is an independent copy of $G(n, p')$
with $(1 - p')^k = 1 - p$, and find the $k$ perfect matchings sequentially, namely, for each $\ell \in [k]$, we find in $G_\ell$ a perfect matching $M_\ell$ that forms a Hamilton cycle with the previously found matchings $M_1,\dots,M_{\ell-1}$.

Let $1 \le \ell \le k$.
Suppose we have found $\ell-1$ edge-disjoint perfect matchings $M_1, \dots, M_{\ell-1}$ in $G_1 \cup \dots \cup G_{\ell-1}$. 
(We know in addition that the union of any two of them induces a Hamilton cycle, but this is not needed for the following argument.)
Let $\cM:=\{M_1, \dots, M_{\ell-1}\}$.
Define the auxiliary hypergraph $\cH$ with $V(\cH)=E(K_n)$, and where each hyperedge of $\cH$ corresponds to a red perfect matching of $\cM$.
Our goal is to show that $\cH$ is $q$-spread for $q:=\gamma n$ for some constant $\gamma>0$ (which can possibly depend on $\ell$). Then, by Theorem~\ref{thm:spread}, w.h.p.~$G_\ell$ contains a matching with the desired property.

By \cref{thm:count robust}, we know that there exist constants $\alpha,\beta$ (depending on $\ell$) with $\alpha \le \beta$ such that for any $(\ell-1)$-configuration on $n$ vertices with total overlap at most $n^{1/8}$, the number of perfect matchings which form a Hamilton cycle with each matching of the configuration is at least $\alpha \cdot n^{-(\ell-1)/2}(n/e)^{n/2}$ and at most $\beta \cdot n^{-(\ell-1)/2}(n/e)^{n/2}$.
We show that choosing $\gamma:=\alpha/(3\beta\cdot e)$ suffices.

Showing that $\cH$ is $q$-spread requires to prove that for every $S \subseteq V(\cH)$, 
\begin{enumerate}[label = (\ding{61})]
	\item \label{item:spread} the number of hyperedges containing $S$ is at most $q^{-|S|}\cdot |\cH|$.
\end{enumerate}
Note that $S$ is a set of edges of~$K_n$ and, if these edges do not form a partial red matching of $\cM$, then no hyperedge of $\cH$ contains~$S$ and there is nothing to prove. So we may assume that $S$ is a partial red matching, and we let $\cM'$ be the reduced configuration of $S$.
Then the number of hyperedges containing $S$ is equal to the number of red perfect matchings of $\cM'$.

By Theorem~\ref{thm:count robust}, we know that $\cH$ has at least $\alpha \cdot n^{-(\ell-1)/2}(n/e)^{n/2}$ hyperedges.
Let $s:=|S|$ and set $n':=n-2s$. 
If $s\ge k \log n$, it is enough to upper bound the number of red perfect matchings of $\cM'$ by the number of all perfect matchings on a set of size $n'$, which in turn is $O((n'/e)^{n'/2}) = O((n/e)^{n/2-s})$.
In order to prove~\ref{item:spread} for $S$, it is enough to show that
\begin{align*}
    O((n/e)^{n/2-s}) \le q^{-s} \cdot \alpha n^{-(\ell-1)/2}(n/e)^{n/2}\, .
\end{align*}  
By plugging in the value of $q$, this is equivalent to $O(n^{(\ell-1)/2}) \le \alpha \cdot \left( \frac{3\beta}{\alpha}\right)^s$, which holds as $\left( \frac{3\beta}{\alpha}\right)^s \ge 3^s \ge e^s \ge e^{k \log n} = n^k = \omega(n^{(\ell-1)/2})$.

Assume now that $s\le k\log n$.
By Item~\ref{prop_P_2} of \Cref{prop:counting_double_edges}, each time we add a red edge, the total overlap can increase by at most $2 \binom{k}{2} \le k^2$.
Since the matchings of $\cM$ are pairwise edge-disjoint, we have $R(\cM)=0$ and thus $R(\cM') \le R(\cM)+k^2s=k^2s=O_k(\log n)$. 
Moreover $\cM'$ has order $n' \ge n/2$.
Therefore, again by Theorem~\ref{thm:count robust}, the number of red perfect matchings of $\cM'$ is at most
\begin{align*}
\beta \cdot n'^{-(\ell-1)/2}(n'/e)^{n'/2} \le \beta \cdot n'^{-(\ell-1)/2}(n/e)^{n/2-s}\, .
\end{align*}
In order to prove~\ref{item:spread} for $S$, it is enough to show that
\begin{align*}
\beta \cdot n'^{-(\ell-1)/2}(n/e)^{n/2-s}  \le q^{-s} \cdot \alpha n^{-(\ell-1)/2}(n/e)^{n/2}.
\end{align*}
By plugging in the value of $q$, this is equivalent to $\left( \frac{n}{n'}\right)^{(\ell-1)/2} \le 3^s \cdot \left(\beta/\alpha\right)^{s-1}$, which can be shown to hold as follows.
We have $\frac{n}{n'}= \frac{1}{1-2s/n} \le 1+\frac{3s}{n} \le \exp\left(\frac{3s}{n}\right)$ and thus $\left( \frac{n}{n'}\right)^{(\ell-1)/2} \le \exp\left(\frac{3(\ell-1)}{2n} \cdot s\right) \le 3^s$ for $n$ large enough.
Moreover, $\beta \ge \alpha$ and thus $3^s \le 3^s \cdot \left(\beta/\alpha\right)^{s-1}$. 

This proves~\ref{item:spread} for all $S$.
\endproof

\section{Bipartite setting}
\label{sec:bipartite}
In light of the conjecture of Wanless~\cite{wanless:99} on a bipartite version of Kotzig's conjecture (see \cref{sec:intro}), it is natural to ask if our results can be adapted to the bipartite setting.
Let $k \in \mathbb{N}$ and $M_1, \dots, M_k$ be perfect matchings of $K_{n,n}$ with small overlap.
Is there a perfect matching $M^*$ (of $K_{n,n}$) which creates a Hamilton cycle with each of them?
And, if so, how many?

A perfect matching $M$ of $K_{n,n}$ can be expressed as the permutation $\pi_M$ of $S_n$, where $\pi_M(i)=j$ if and only if $ij \in M$.
Observe that if there exists $M^*$ such that $M^* \cup M_i$ is a Hamilton cycle, then $\pi_{M^*} \circ \pi_{M_i}$ consists of one cycle and its parity depends on $n$ only (it is even if $n$ is odd and odd otherwise).
Since this is true for each $i \in [k]$, it follows that the $\pi_{M_i}$ must all have the same parity.
In particular, for all distinct $i,j \in [k]$, the permutation $\pi_{M_i} \circ \pi_{M_j}$ is even and thus 
\begin{enumerate}[label = (\ding{69})]
	\item \label{item:bipartite} $M_i \cup M_j$ has an even number of cycles whose length is divisible by $4$.
\end{enumerate}
This shows that condition~\ref{item:bipartite} is necessary for the existence of such $M^*$. 
Ironically, among our arguments, the only one which requires a conceptual modification to show that condition~\ref{item:bipartite} is also sufficient is \cref{prop:two colours}, which was the easiest part in the general setting.
Therefore we provide the full details here.
\begin{prop}[Analogue of \cref{prop:two colours} in the bipartite setting]
\label{prop:two colours bipartite}
    Fix two matchings of $K_{n,n}$, say of colour blue and green, such that condition~\ref{item:bipartite} holds.
    Then there exists an equalizing partial red matching.
\end{prop}
In the non-bipartite setting, this was achieved by inserting red edges inside the alternating blue-green cycles to break them into shorter alternating cycles and double edges. In particular, when the cycle had length exactly $4$, we were taking one of the diagonals as a red edge.
This is no longer possible in the bipartite setting as the diagonals are non-bipartite edges.
However, we own the extra assumption~\ref{item:bipartite} and use the following fact.

\begin{fact}
\label{fact:bipartite}
    Given two perfect matchings of $K_{n,n}$, consider the cycles in their union.
    Then the parity of the number of those whose length is divisible by $4$ is an invariant under the addition of (bipartite) red edges.
\end{fact}

\lateproof{\cref{prop:two colours bipartite}}
The union of the two matchings consists of alternating blue-green even cycles and double edges. For the latter case, we do not have to do anything, so suppose we have a cycle with edges alternating blue and green. If the length is at least $6$, then taking a red edge between two vertices at distance $3$ on the cycle splits off a double-coloured edge and an alternating cycle whose length is shorter by $4$.
(In particular, a cycle of length $6$ splits off into two double edges.)
By repeatedly choosing red edges in this way, we get a reduced configuration consisting of alternating blue-green cycles of length $4$ and double edges.
By \cref{fact:bipartite}, condition~\ref{item:bipartite} still holds and thus there is an even number of cycles of length $4$.
We pair these cycles arbitrarily and add a (bipartite) red edge between each pair. 
This leaves a reduced configuration consisting of alternating cycles of length $6$ and double edges, which can be easily equalized.
\endproof

The rest of the arguments goes through and we can prove the following results for $K_{n,n}$.

\begin{theorem}
For every fixed $k\in \bN$, there exists a constant $C=C(k)$ such that the following holds. Let $M_1,\dots,M_k$ be perfect matchings of $K_{n,n}$ satisfying condition~\ref{item:bipartite} and such that
\begin{align*}
|M_i\sm (\cup_{j\neq i}M_j)| \ge C
\end{align*}
for each $i\in[k]$.
Then there exists a perfect matching $M^*$ of $K_{n,n}$ such that $M^*\cup M_i$ forms a Hamilton cycle for each $i\in [k]$.
\end{theorem}

\begin{theorem}
\label{thm:counting bipartite}
    Let $k\in \bN$ be fixed.
    Then the following holds for all $n$ sufficiently large. Let $M_1,\dots,M_k$ be perfect matchings of $K_{n,n}$ satisfying condition~\ref{item:bipartite} and such that $\sum_{i,j:i\neq j} |M_i\cap M_j| \le n^{1/8}$.
    Then the number of perfect matchings $M^*$ of $K_{n,n}$ such that $M^*\cup M_i$ forms a Hamilton cycle for each $i\in [k]$ is $\Theta_k(n^{1/2-k} (n/e)^n)$. 
\end{theorem}

Finally, we discuss the analogue of \cref{thm:random} for the bipartite random graph $G(n,n,p)$.
For a fixed $k \in \mathbb{N}$, under which condition on $p$ does $G(n,n,p)$ contain $k$ perfect matchings such that any pair forms a Hamilton cycle?
If $k=2$, the question is already well understood as we are asking for a Hamilton cycle.
So we can assume $k \ge 3$.
First observe that if $M$ and $M'$ are perfect matchings of $K_{n,n}$ and $M \cup M'$ is a Hamilton cycle, then, as above, $\pi_M \circ \pi_{M'}$ is an odd permutation if $n$ is even (and it is even otherwise).
This implies that, if $n$ is even, then $\pi_M$ and $\pi_{M'}$ have different parity.
In particular, if $K_{n,n}$ admits three (or more) perfect matchings such that any pair forms a Hamilton cycle, then $n$ must be odd.
Therefore, for our question to have a positive answer, this condition is necessary.
We can prove the following.
\begin{theorem}
    For every fixed $k \in \mathbb{N}$ with $k \ge 3$, there exists a constant $C=C(k)$ such that, when $p\ge \frac{C \log n}{n}$, w.h.p.~for $n$ odd the bipartite random graph $G(n,n,p)$ contains $k$ edge-disjoint perfect matchings $M_1, \dots, M_k$ such that $M_i \cup M_j$ induces a Hamilton cycle for all distinct $i,j \in [k]$.
\end{theorem}
The proof follows the steps of that of \cref{thm:random} and, in particular, builds the matchings iteratively one after the other in $k$ rounds, using \cref{thm:counting bipartite} and \cref{thm:spread}.
Note that if $n$ is odd and a collection of perfect matchings of $K_{n,n}$ is such that any pair forms a Hamilton cycle, then condition~\ref{item:bipartite} clearly holds and thus \cref{thm:counting bipartite} can be applied.

\section{Concluding remarks}\label{sec:remarks}
\cref{thm:random} offers a variant of Kotzig's conjecture for random graphs.
A natural problem is to sharpen the threshold for the edge probability. One might expect that $p\ge \frac{\log n+(k-1)\log\log n + \omega(n)}{n}$ is already enough.
In fact, even a ``hitting time'' result might be possible. The legendary hitting time result, proved independently by Bollob\'as~\cite{bollobas:84} and Ajtai, Koml\'os and Szemer\'edi~\cite{AKS:85}, says that if we start with the empty graph and then keep adding edges at random, w.h.p., in the very moment that the resulting graph has minimum degree $2$, it will also have a Hamilton cycle.
Later, Bollob\'as and Frieze~\cite{BF:85} extended this result by proving that, for a given constant $k$, in the very moment that the minimum degree is $2k$, there exist $k$ edge-disjoint Hamilton cycles. The following would constitute a powerful generalization.

\begin{problem}
Prove that, in the random graph process, w.h.p., as soon as the minimum degree is $k$, there exists a collection of $k$ edge-disjoint perfect matchings such that the union of any two distinct of them induces a Hamilton cycle.
\end{problem}
Note that this result would imply the result of Bollob\'as and Frieze~\cite{BF:85} since we can arbitrarily pair up the $2k$ perfect matchings and each pair will form a Hamilton cycle, and the $k$ resulting Hamilton cycles will be edge-disjoint.
A closely related problem would be to consider random regular graphs. We reiterate a conjecture by Wormald~\cite{wormald:99} that for fixed $d \ge 3$, a random $d$-regular graph should w.h.p.~have a perfect $1$-factorisation.

\medskip

Another direction for future work is to sharpen the counting result (\cref{thm:count}) to obtain an asymptotic estimate.
Given the heuristic argument described in \cref{sec:intro}, we conjecture that given $k$ edge-disjoint perfect matchings on $[n]$, the number of perfect matchings forming a Hamilton cycle with each of them is $\left(\sqrt{2} \cdot \left(\pi/2\right)^{k/2} \pm o(1)\right) \cdot n^{-k/2} \cdot (n/e)^{n/2}$.
We recall that the case $k=2$ has already been proved by Kim and Wormald~\cite{KW:01}.
We remark that, even with a more precise counting result, our method would not give a sharper result for $G(n,p)$ since it relies on the expectation threshold conjecture.

\section*{Acknowledgments}
We are grateful to Pranshu Gupta for stimulating discussions at the early stages of this project.
We would like to thank the anonymous referees for a careful reading of the paper and one of them for pointing us to an inaccuracy in the proof of Lemma 4.5 in an
earlier version of this manuscript.
\bibliographystyle{amsplain_v2.0customized}
\bibliography{References}

\end{document}